\documentclass[12pt]{article}
\usepackage{latexsym, amssymb}

\DeclareMathAlphabet\gothic{U}{euf}{m}{n}

\makeatletter
\def\eqnarray{\stepcounter{equation}\let\@currentlabel=\theequation
\global\@eqnswtrue
\tabskip\@centering\let\\=\@eqncr
$$\halign to \displaywidth\bgroup\hfil\global\@eqcnt\z@
  $\displaystyle\tabskip\z@{##}$&\global\@eqcnt\@ne
  \hfil$\displaystyle{{}##{}}$\hfil
  &\global\@eqcnt\tw@ $\displaystyle{##}$\hfil
  \tabskip\@centering&\llap{##}\tabskip\z@\cr}

\def\endeqnarray{\@@eqncr\egroup
      \global\advance\c@equation\m@ne$$\global\@ignoretrue}

\def\@yeqncr{\@ifnextchar [{\@xeqncr}{\@xeqncr[5pt]}}
\makeatother

\begin{document}
\bibliographystyle{tom}

\newtheorem{lemma}{Lemma}[section]
\newtheorem{thm}[lemma]{Theorem}
\newtheorem{cor}[lemma]{Corollary}
\newtheorem{voorb}[lemma]{Example}
\newtheorem{rem}[lemma]{Remark}
\newtheorem{prop}[lemma]{Proposition}
\newtheorem{stat}[lemma]{{\hspace{-5pt}}}
\newtheorem{obs}[lemma]{Observation}
\newtheorem{defin}[lemma]{Definition}

\newenvironment{remarkn}{\begin{rem} \rm}{\end{rem}}
\newenvironment{exam}{\begin{voorb} \rm}{\end{voorb}}
\newenvironment{defn}{\begin{defin} \rm}{\end{defin}}
\newenvironment{obsn}{\begin{obs} \rm}{\end{obs}}

\newenvironment{emphit}{\begin{itemize} }{\end{itemize}}

\newcommand{\gota}{\gothic{a}}
\newcommand{\gotb}{\gothic{b}}
\newcommand{\gotc}{\gothic{c}}
\newcommand{\gote}{\gothic{e}}
\newcommand{\gotf}{\gothic{f}}
\newcommand{\gotg}{\gothic{g}}
\newcommand{\gothh}{\gothic{h}}
\newcommand{\gotk}{\gothic{k}}
\newcommand{\gotm}{\gothic{m}}
\newcommand{\gotn}{\gothic{n}}
\newcommand{\gotp}{\gothic{p}}
\newcommand{\gotq}{\gothic{q}}
\newcommand{\gotr}{\gothic{r}}
\newcommand{\gots}{\gothic{s}}
\newcommand{\gotu}{\gothic{u}}
\newcommand{\gotv}{\gothic{v}}
\newcommand{\gotw}{\gothic{w}}
\newcommand{\gotz}{\gothic{z}}
\newcommand{\gotA}{\gothic{A}}
\newcommand{\gotB}{\gothic{B}}
\newcommand{\gotG}{\gothic{G}}
\newcommand{\gotL}{\gothic{L}}
\newcommand{\gotS}{\gothic{S}}
\newcommand{\gotT}{\gothic{T}}

\newcommand{\mn}{\marginpar{\hspace{1cm}*} }
\newcommand{\mnn}{\marginpar{\hspace{1cm}**} }

\newcommand{\mnq}{\marginpar{\hspace{1cm}*???} }
\newcommand{\mnnq}{\marginpar{\hspace{1cm}**???} }

\newcounter{teller}
\renewcommand{\theteller}{\Roman{teller}}
\newenvironment{tabel}{\begin{list}%
{\rm \bf \Roman{teller}.\hfill}{\usecounter{teller} \leftmargin=1.1cm
\labelwidth=1.1cm \labelsep=0cm \parsep=0cm}
                      }{\end{list}}

\newcounter{tellerr}
\renewcommand{\thetellerr}{(\roman{tellerr})}
\newenvironment{subtabel}{\begin{list}%
{\rm  (\roman{tellerr})\hfill}{\usecounter{tellerr} \leftmargin=1.1cm
\labelwidth=1.1cm \labelsep=0cm \parsep=0cm}
                         }{\end{list}}
\newenvironment{ssubtabel}{\begin{list}%
{\rm  (\roman{tellerr})\hfill}{\usecounter{tellerr} \leftmargin=1.1cm
\labelwidth=1.1cm \labelsep=0cm \parsep=0cm \topsep=1.5mm}
                         }{\end{list}}

\newcommand{\Ni}{{\bf N}}
\newcommand{\Ri}{{\bf R}}
\newcommand{\Ci}{{\bf C}}
\newcommand{\Ti}{{\bf T}}
\newcommand{\Zi}{{\bf Z}}
\newcommand{\Fi}{{\bf F}}

\newcommand{\proof}{\mbox{\bf Proof} \hspace{5pt}} 
\newcommand{\remark}{\mbox{\bf Remark} \hspace{5pt}}
\newcommand{\ruimte}{\vskip10.0pt plus 4.0pt minus 6.0pt}

\newcommand{\simh}{{\stackrel{{\rm cap}}{\sim}}}
\newcommand{\ad}{{\mathop{\rm ad}}}
\newcommand{\Ad}{{\mathop{\rm Ad}}}
\newcommand{\Aut}{\mathop{\rm Aut}}
\newcommand{\arccot}{\mathop{\rm arccot}}
\newcommand{\capp}{{\mathop{\rm cap}}}
\newcommand{\rcapp}{{\mathop{\rm rcap}}}
\newcommand{\Capp}{{\mathop{\rm Cap}}}
\newcommand{\diam}{\mathop{\rm diam}}
\newcommand{\divv}{\mathop{\rm div}}
\newcommand{\dist}{\mathop{\rm dist}}
\newcommand{\codim}{\mathop{\rm codim}}
\newcommand{\RRe}{\mathop{\rm Re}}
\newcommand{\IIm}{\mathop{\rm Im}}
\newcommand{\Tr}{{\mathop{\rm Tr}}}
\newcommand{\Vol}{{\mathop{\rm Vol}}}
\newcommand{\card}{{\mathop{\rm card}}}
\newcommand{\supp}{\mathop{\rm supp}}
\newcommand{\sgn}{\mathop{\rm sgn}}
\newcommand{\essinf}{\mathop{\rm ess\,inf}}
\newcommand{\esssup}{\mathop{\rm ess\,sup}}
\newcommand{\Int}{\mathop{\rm Int}}
\newcommand{\Leibniz}{\mathop{\rm Leibniz}}
\newcommand{\lcm}{\mathop{\rm lcm}}
\newcommand{\loc}{{\rm loc}}

\newcommand{\mod}{\mathop{\rm mod}}
\newcommand{\spann}{\mathop{\rm span}}
\newcommand{\one}{1\hspace{-4.5pt}1}

\newcommand{\DWR}{}

\hyphenation{groups}
\hyphenation{unitary}

\newcommand{\tfrac}[2]{{\textstyle \frac{#1}{#2}}}

\newcommand{\cb}{{\cal B}}
\newcommand{\cc}{{\cal C}}
\newcommand{\cd}{{\cal D}}
\newcommand{\ce}{{\cal E}}
\newcommand{\cf}{{\cal F}}
\newcommand{\ch}{{\cal H}}
\newcommand{\ci}{{\cal I}}
\newcommand{\ck}{{\cal K}}
\newcommand{\cl}{{\cal L}}
\newcommand{\cm}{{\cal M}}
\newcommand{\cn}{{\cal N}}
\newcommand{\co}{{\cal O}}
\newcommand{\cs}{{\cal S}}
\newcommand{\ct}{{\cal T}}
\newcommand{\cx}{{\cal X}}
\newcommand{\cy}{{\cal Y}}
\newcommand{\cz}{{\cal Z}}

\newcommand{\wtozp}{W^{1,2}\raisebox{10pt}[0pt][0pt]{\makebox[0pt]{\hspace{-34pt}$\scriptstyle\circ$}}}
\newlength{\hightcharacter}
\newlength{\widthcharacter}
\newcommand{\covsup}[1]{\settowidth{\widthcharacter}{$#1$}\addtolength{\widthcharacter}{-0.15em}\settoheight{\hightcharacter}{$#1$}\addtolength{\hightcharacter}{0.1ex}#1\raisebox{\hightcharacter}[0pt][0pt]{\makebox[0pt]{\hspace{-\widthcharacter}$\scriptstyle\circ$}}}
\newcommand{\cov}[1]{\settowidth{\widthcharacter}{$#1$}\addtolength{\widthcharacter}{-0.15em}\settoheight{\hightcharacter}{$#1$}\addtolength{\hightcharacter}{0.1ex}#1\raisebox{\hightcharacter}{\makebox[0pt]{\hspace{-\widthcharacter}$\scriptstyle\circ$}}}
\newcommand{\scov}[1]{\settowidth{\widthcharacter}{$#1$}\addtolength{\widthcharacter}{-0.15em}\settoheight{\hightcharacter}{$#1$}\addtolength{\hightcharacter}{0.1ex}#1\raisebox{0.7\hightcharacter}{\makebox[0pt]{\hspace{-\widthcharacter}$\scriptstyle\circ$}}}

\newpage

 \thispagestyle{empty}
 
 \begin{center}

 \vspace*{-1.0cm}
\vspace*{1.5cm}
 
{\Large{\bf Hardy and Rellich inequalities}  }\\[3mm]
{\Large{\bf on the complement of convex sets }}  \\[4mm]
\large Derek W. Robinson$^\dag$ \\[1mm]

\normalsize{April 2017}
\end{center}

\vspace{+5mm}

\begin{list}{}{\leftmargin=1.7cm \rightmargin=1.7cm \listparindent=15mm 
   \parsep=0pt}
   \item
{\bf Abstract} $\;$ 
We establish existence of  weighted Hardy and Rellich inequalities on  the spaces $L_p(\Omega)$ where  $\Omega= \Ri^d\backslash K$ with $K$  a closed convex subset of $\Ri^d$.
Let $\Gamma=\partial\Omega$ denote the boundary of $\Omega$ and $d_\Gamma$ the Euclidean distance to $\Gamma$.
We consider weighting functions $c_\Omega=c\circ d_\Gamma$ with  $c(s)=s^\delta(1+s)^{\delta'-\delta}$ and $\delta,\delta'\geq0$.
Then the  Hardy inequalities take the form
\[
\int_\Omega c_\Omega\,|\nabla\varphi|^p\geq   b_p\int_\Omega c_\Omega\,d_\Gamma^{\;-p}\,|\varphi|^p
\]
and the Rellich inequalities are given by
\[
\int_\Omega|H\varphi|^p\geq d_p\int_\Omega |c_\Omega\,d_\Gamma^{\,-2}\varphi|^p
\]
with $H=-\divv(c_\Omega\nabla)$.
The constants $b_p, d_p$ depend on the weighting parameter $\delta,\delta'\geq0$ and the Hausdorff dimension of the boundary. 
We compute the optimal constants in a broad range of situations.

\end{list}

\vfill


\vspace{0.5cm}

\noindent
\begin{tabular}{@{}cl@{\hspace{10mm}}cl}
$ {}^\dag\hspace{-5mm}$&   Mathematical Sciences Institute (CMA)    &  {} &{}\\
  &Australian National University& & {}\\
&Canberra, ACT 0200 && {} \\
  & Australia && {} \\
  &derek.robinson@anu.edu.au
 & &{}\\
\end{tabular}

\newpage

\setcounter{page}{1}

\section{Introduction}\label{S1}

The classical Hardy and Rellich inequalities are  estimates for differential operators on the spaces $L_p(\Ri^d\backslash\{0\})$, $p\in\langle1,\infty\rangle$,  for which the optimal constants are known.
Our intention is to derive similar  estimates on $\Omega=\Ri^d\backslash K$ where $K$ is a closed convex subset of $\Ri^d$.
There are two stages in the analysis, first the existence of the inequalities and secondly the optimality of the corresponding constants.
Background information on both these aspects and references to the literature can be found in the recent monograph \cite{BEL}.
Our  estimates depend on the  Hausdorff--Minkowski dimension $d_H$ of the boundary $\Gamma=\partial\Omega$ of $\Omega$. 
If the dimension $\dim (K)$ (of the affine closure $A_K$) of $K$ 
 takes one of the values $0,1,\ldots,d-1$ then $ d_H=\dim(K)$ but if $\dim(K)=d$ then $d_H=d-1$, assuming
 $K\neq\Ri^d$.
In addition, the general inequalities depend on
the Euclidean distance $d_\Gamma$ to the boundary, i.e.\
$ d_\Gamma(x)=\inf_{y\in \Omega^c}|x-y|$ for $x\in \Omega$.
 We begin by  establishing the existence of  weighted Hardy inequalities with a weight function $c_\Omega=c\circ d_\Gamma$ where   $c$ is  a strictly positive function  on $\langle0,\infty\rangle$ with different power behaviours at the origin and at infinity.

\begin{thm}\label{tcc1.1}
Let $\Omega=\Ri^d\backslash K$ where $K$ is a closed  convex subset of $\Ri^d$
and  denote  the Hausdorff dimension of the boundary $\Gamma$ of $ \Omega$ by $d_H$.
Further let $c_\Omega=c\circ d_\Gamma$ where  $c(s)=s^\delta(1+s)^{\delta-\delta'}$ with $\delta,\delta'\geq0$.
If $d-d_H+(\delta\wedge\delta')-p>0$ with $p\in[1,\infty\rangle$ then
\begin{equation}
\int_\Omega c_\Omega\,|\nabla\varphi|^p\geq 
\int_\Omega c_\Omega\,|(\nabla d_\Gamma).(\nabla\varphi)|^p\geq a_p^{\,p}\int_\Omega c_\Omega\,d_\Gamma^{\;-p}\,|\varphi|^p
\label{ecc1.1}
\end{equation}
for all $\varphi\in C_c^1(\Omega)$ with $a_p=(d-d_H+(\delta\wedge\delta')-p)/p$.
\end{thm}

Here and in the sequel all functions are real-valued. 
Moreover, we use the standard notation $|\nabla\varphi|=(\sum^d_{k=1}|\partial_k\varphi|^2)^{1/2}$.
Then the left hand inequality in (\ref{ecc1.1}) follows since
$d_\Gamma$ is a Lipschitz function with $|\nabla d_\Gamma|\leq1$.
The choice of the weight $c$ is governed by the asymptotic properties $c(s)/s^\delta\to1$ as $s\to0$ and $c(s)/s^{\delta'}\to 1$ as $s\to\infty$.
Although this  theorem and the subsequent one are stated for the particular coefficient $c$ the general conclusions are valid for a large 
class of $c$ with similar asymptotic properties.
Note that if $\delta=\delta'$ then $c(s)=s^\delta$ which is the conventional weight function  used in the discussion of Hardy inequalities.

An important part of the proof of the theorem, which will be given in Section~\ref{S2}, is the observation that $d_\Gamma$ is a convex function on convex subsets of~$\Omega$.
This is the analogue of the statement that if $U$ is an open  convex subset of $\Ri^d$   then $d_{\partial U}$ is  a concave function (see
\cite{Hor8}, Corollary~2.1.26, or  \cite{BFT}, Example~2).
The proofs of the two statements are very  similar.
There is also a somewhat weaker  analogue of Theorem~\ref{tcc1.1} for weighted operators on convex sets which we will establish at the end of Section~\ref{S2}.

In Section~\ref{S3} we consider the existence of weighted Rellich inequalities on $\Omega=\Ri^d\backslash K$.
The classic Rellich inequalities were initially established for the Laplacian $\Delta=-\sum^d_{k=1}\partial_k^{\,2}$ on $L_2(\Ri^d\backslash\{0\})$ but have subsequently been extended to all the spaces $L_p(\Ri^d\backslash\{0\})$ with $p\in\langle1,\infty\rangle$
(see, for example, \cite{BEL} Sections~6.1--6.3 and in particular Corollary~6.3.5).
Our aim is to establish similar estimates for the weighted operators $H=-\sum^d_{k=1}\partial_k \,c_\Omega\,\partial_k=-\divv(c_\Omega\nabla)$ on the spaces $L_p(\Omega)$.
The operators $H$ are  defined on the universal domain $C_c^2(\Omega)$ and  all estimates are on this domain.

\begin{thm}\label{tcc1.2}
Let $\Omega=\Ri^d\backslash K$ where $K$ is a closed  convex subset of $\Ri^d$
and  denote  the Hausdorff dimension of the boundary $\Gamma$ of $\Omega$ by  $d_H$.
Further let $c_\Omega=c\circ d_\Gamma$ where  $c(s)=s^\delta(1+s)^{\delta-\delta'}$ with $\delta,\delta'\in[0,2\rangle$.
If $d-d_H+p\,(\delta\wedge\delta')-2p\geq2p\,|\delta-\delta'|\,(2-\delta\vee\delta')^{-1}$
with $p\in\langle1,\infty\rangle$ then there is a $c_p\in\langle0,C_p]$, where $C_p=(p-1)\,(d-d_H)\,(d-d_H+p\,(\delta\wedge\delta')-2p)\,p^{-2}$,
such that 
\begin{equation}
\int_\Omega|H\varphi|^p\geq c_p^{\,p}\int_\Omega |c_\Omega\,d_\Gamma^{\,-2}\varphi|^p
\label{ecc1.2}
\end{equation}
for all $\varphi\in C_c^2(\Omega)$.
Moreover, if $\delta=\delta'$ then $c_p=C_p$. 
\end{thm}

The proof of the theorem  allows one to identify $c_p$ as a function of $d_H, \delta$ and $\delta'$
but the result is significantly more complicated than the expression for $C_p$.
Although the condition $c_p>0$ requires the restriction  $\delta,\delta'<2$   the existence of the Rellich inequalities should not depend on this latter condition.
In fact if $p=2$ then the weighted inequalities (\ref{ecc1.2}) follow for  all $\delta,\delta'\geq 0$ from the arguments of  \cite{Rob12}.

Theorems~\ref{tcc1.1} and \ref{tcc1.2} establish criteria for existence of the Hardy and Rellich inequalities (\ref{ecc1.1}) and 
(\ref{ecc1.2}), respectively, but they give no information about optimality of the constants $a_p^{\,p}$ and $c_p^{\,p}$.
This problem, which appears more challenging, is tackled in Section~\ref{S4}.
We show, for example, that $a_p^{\,p}$ is optimal for the Hardy inequality if $K=\{0\}$.
It is also optimal   if $k=\dim(K)\in\{1,\ldots,d-1\}$
and either $\delta\leq \delta'$ or  the `dimension at infinity'  $k_\infty$ of $K$ is equal to $k$.
Alternatively $C_p^{\,p}$ is the optimal constant for the Rellich inequality if $K=\{0\}$ and  $\delta=\delta'\in[0,2\rangle$.
More generally it is optimal if $k\in\{1,\ldots,d-1\}$, $\delta=\delta'\in[0,2\rangle$ and  $k_\infty=k$.
But these results leave open room for improvement.
In particular if  $p=2$ then  it follows from  \cite{Rob12} that  $C_2^{\,2}$ is the optimal constant for all $\delta,\delta'\geq0$ 
such that $\delta+\delta'\leq4$ with no further restriction on $K$ or $\delta,\delta'$ other than $C_2>0$.

Finally note that if  there is no weighting factor, i.e.\ if $c=1$, then $H=\Delta$, the Laplacian, and (\ref{ecc1.2}) states that
\[
\int_\Omega|\Delta\varphi|^p\geq C_p^{\,p}\int_\Omega |d_\Gamma^{\,-2}\varphi|^p
\]
for all $\varphi\in C_c^2(\Omega)$.
Moreover, the constant $C_p=(p-1)(d-d_H)(d-d_H-2p)/p^2$ is optimal.
In particular if $K=\{0\}$ then $d_H=0$ and this gives the classical $L_p$-Rellich inequality with the optimal constant.
But (\ref{ecc1.2}), which is a `weighted operator' version of the classical inequality,  is not the only possible weighted generalization.
A second natural alternative would be the `weighted measure' version
\begin{equation}
\int_\Omega c_\Omega^{\;p}\,|\Delta\varphi|^p\geq b_p\int_\Omega c_\Omega^{\;p}\,|d_\Gamma^{\,-2}\varphi|^p
\label{ecc1.3}
\end{equation}
with $b_p>0$ for all $\varphi\in C_c^2(\Omega)$.
The relation  between the existence and optimality  of the two versions (\ref{ecc1.2}) and (\ref{ecc1.3}) of the Rellich inequalities is unclear.

\section{Hardy inequalities}\label{S2}

In this section we prove Theorem~\ref{tcc1.1}.
As a preliminary to the proof we need to establish local convexity of the distance function $d_\Gamma$
where  $\Omega=\Ri^d\backslash K$ with $K$ a closed convex subset $K$  of $\Ri^d$ and  $K\neq\Ri^d$.
Since  $K$ is the complement of $\Omega$
it follows from Motzkin's theorem (see, for example, \cite{Hor8},  Theorem~2.1.20, or \cite{BEL}, Theorem~2.2.9) that each point $x\in \Omega$ has  a unique nearest point $n(x)\in K$, i.e.\ there is a unique $n(x)\in K$ such that  $d_\Gamma(x)=|x-n(x)|$.
Moreover, $d_\Gamma$ is differentiable at each point $x\in\Omega$ and $(\nabla d_\Gamma)(x)=(x-n(x))/|x-n(x)|$.
 Thus $|\nabla d_\Gamma|=1$ and  $(\nabla d_\Gamma^{\;2})(x)=2\,(x-n(x))$.
Note that in  the degenerate case  $K=\{0\}$ one has  $d_\Gamma(x)=|x|$
and consequently  $\nabla^2d_\Gamma^{\;2}=2\,d$.
In the non-degenerate case it  is not, however, clear that $d_\Gamma$ is even twice-differentiable.
But this follows from local convexity.
 
 \begin{prop}\label{pcc2.1}
 The distance function $d_\Gamma$ is convex on all open convex subsets of $\Omega$.
 In particular it is twice-differentiable almost everywhere in $\Omega$ and the corresponding Hessian
 $(\partial_k\partial_ld_\Gamma)(x)$ is positive-definite for 
almost all $x\in\Omega$.
 \end{prop}
 \proof\
 First we prove the convexity in an open neighbourhood of an arbitrarily chosen point
of  $\Omega$.

 Let   $n(x)\in \Gamma$ be the unique near point of  $x\in\Omega$.
Then there is a unique tangent hyperplane $T_x$ at the point $n(x)$ which is orthogonal to $x-n(x)$.
 The hyperplane separates $\Ri^d$ into two half open half spaces, $\Gamma^{(+)}_x\subset \Omega$ and 
 $\Gamma^{(-)}_x\supset \Int(\Omega^c)$.
 Moreover, $\Omega=\bigcup_{x\in\Omega}\Gamma^{(+)}_x$ and $\Int(\Omega^c)=\bigcap_{x\in \Omega}\Gamma^{(-)}_x$.
 Now fix a point $x_0\in\Omega$ and an $r>0$ such that the open  Euclidean ball $B_{x_0}(r)$ with centre $x_0$ and radius $r$
 is contained in $\Omega$.
 Next choose $r$ sufficiently small that $B_{x_0}(r)\subset \bigcap_{x\in B_{x_0}(r)}\Gamma^{(+)}_x$.
 This is possible since if $x_k\in \Omega$ converges pointwise to $x\in\Omega$ then $n(x_k)\to n(x)$ (see \cite{BEL}, Lemma~2.2.1).
 Therefore the family of open subsets $s>0\mapsto \Lambda_{x_0}(s)= \bigcap_{x\in B_{x_0}(s)}\Gamma^{(+)}_x$
 increases  as $s$ decreases to zero to $\Gamma^{(+)}_{x_0}\supset B_{x_0}(r)$.
 But the balls $B_{x_0}(s)$ decrease as $s\to0$.
 Therefore there is an $r_0$ such that $B_{x_0}(r)\subset \bigcap_{x\in B_{x_0}(r_0)}\Gamma^{(+)}_x$ for all $r\in\langle0,r_0\rangle$.
 
 Secondly, we argue that if $r<r_0$ then $d_\Gamma$ is convex on $B_{x_0}(r)$.
 To this end choose three points $x, y,z\in B_{x_0}(r)$ such that $x=\lambda\,y+(1-\lambda)\,z$
with  $\lambda\in\langle0,1\rangle$.
 Since $r<r_0$ it follows that $B_{x_0}(r)\subset \Gamma^{(+)}_x$.
 Thus the tangent plane $T_x$ separates $B_{x_0}(r)$ and $\Omega^c$.
 Next let $\tilde x,\tilde y,\tilde z$ denote the orthogonal projections of $x,y,z$ onto $T_x$.
 Then $\tilde x=n(x)$, by definition, and $d_\Gamma(x)=|x-\tilde x|$.
 But 
 \[
 |y-\tilde y|=\inf_{y_0\in \Gamma^{(-)}_x}|y-y_0|\leq \inf_{y_0\in \Omega^c}|y-y_0|=d_\Gamma(y)
 \;.
 \]
 Similarly $|z-\tilde z|\leq d_\Gamma(z)$.
 Moreover, $\tilde x=\lambda\,\tilde y +(1-\lambda)\,\tilde z$ and 
 \[
 |x-\tilde x|=\lambda\,|y-\tilde y|+(1-\lambda)\,|z-\tilde z|
 \;.
 \]
 Therefore $d_\Gamma(x)\leq \lambda\,d_\Gamma(y)+(1-\lambda)\,d_\Gamma(z)$.
 Since this is valid for all choices of $x,y,z\in B_{x_0}(r)$ and $\lambda\in\langle0,1\rangle$  with $x=\lambda\, y+(1-\lambda)\,z$ it follows that $d_\Gamma$ is convex on 
 $B_{x_0}(r)$.
 
 Thirdly, it follows from Motzkin's theorem that $d_\Gamma$ is once-differentiable at each $x\in\Omega$.
 But since $d_\Gamma$ is convex on $B_{x_0}(r)$ it follows from Alexandrov's theorem (see \cite{EvG}, Section~6.4) that $d_\Gamma$ is
 twice-differentiable almost-everywhere on $B_{x_0}(r)$.
 Since this is valid for each $x_0\in\Omega$ for some $r>0$
 it then follows that $d_\Gamma$ is twice-differentiable almost-everywhere on $\Omega$.
 The Hessian of a convex function is automatically positive-definite.
 Hence the Hessian of $d_\Gamma$ is positive-definite almost everywhere on $\Omega$.
 
 Finally let $d_\Gamma^{\,(\varepsilon)}$, $\varepsilon>0$, denote a family of local mollifications/regularizations of $d_\Gamma$
 (see \cite{EvG}, Section~4.2.1).
 Then the $d_\Gamma^{\,(\varepsilon)}$ are $C^2$-functions and their Hessians are positive-definite.
 In fact the proof of Alexandrov's theorem relies on  proving the positive-definiteness of the regularizations.
Next  it follows by a standard consequence of convexity (see \cite{bSim7}, Theorem~1.5) that $d_\Gamma^{\,(\varepsilon)}$ is convex on all open convex subsets  suitably distant from the boundary.
 But $d_\Gamma^{\,(\varepsilon)}\to d_\Gamma$ as $\varepsilon\to0$.
 Therefore in the limit $d_\Gamma$ is convex on all open convex subsets of~$\Omega$. \hfill$\Box$
 
 \bigskip
 
 The subsequent proof of the  Hardy inequalities of Theorem~\ref{tcc1.1} depends on control of the second derivatives of $d_\Gamma$.
 \begin{cor}\label{ccc2.2}
  If $\Omega=\Ri^d\backslash K$ where $K$ is a closed convex subset 
   then $\nabla^2d_\Gamma^{\;2}\geq 2\,(d-d_H)$ where $d_H$ is the Hausdorff $($Minkowski$\,)$ dimension of $\Gamma$.
  \end{cor}
\proof\ 
First if $K$ is a singleton then one can assume $K=\{0\}$.
Hence $d_\Gamma^{\;2}(x)=|x|^2$ and $\nabla^2d_\Gamma^{\;2}=2\,d$.

Secondly,   if $\dim(K)=k$ with $k\in \{1,\ldots, d-1\}$  one can factor $\Ri^d$ as a direct product $ \Ri^k\times \Ri^{d-k}$ 
where  $\Ri^k$ is identified with $A_K$, the affine hull of $K$. 
Thus if  $x=(y,z)\in \Ri^d$ with $y\in \Ri^k$ and $z\in \Ri^{d-k}$ one  has $d_\Gamma^{\;2}(x)=d_K^{\;2}(y)+|z|^2$
where $d_K(y)=\inf_{y'\in K}|y-y'|$.
In particular  if $y\in K$ then  $d_\Gamma^{\;2}(x)=|z|^2$
and   $\nabla^2d_\Gamma^{\;2}=\nabla_{\!z}^2d_\Gamma^{\;2}=2\,(d-k)=2\,(d-d_H)$
because $d_H=\dim(K)$.
But if $y\not\in K$ then   $(\nabla_{\!x}^2d_\Gamma^{\;2})(x)=(\nabla_{\!y}^2d_K^{\;2})(y)+\nabla_{\!z}^2|z|^2> 2\,(d-k)$.
Hence one now has  $\nabla^2d_\Gamma^{\;2}\geq 2\,(d-d_H)$ for all  $k\in \{1,\ldots, d-1\}$ .

Thirdly, if $\dim(K)=d$, and $K\neq \Ri^d$ then  $\Gamma=\partial K$ and
the Hausdorff dimension $d_H$ of $\Gamma$ is $d-1$.
Then one can argue as in \cite{BEL}, Section~3.4.2 and 3.4.3 that  $\nabla^2d_\Gamma^{\;2}\geq 2$.
Specifically if $x\in\Omega$ one can choose coordinates $x=(y_1,z)$ with $y_1>0$, $z\in\Ri^{d-1}$
and such that the near point of $(y_1,0)$ is the origin.
Then 
\[
(\nabla_{\!x}^2d_\Gamma^{\;2})(x)=\partial_{y_1}^{\,2}y_1^2+(\nabla_{\!z}^2d_\Gamma^{\;2})(x)\geq2
\]
since $(\nabla_{\!z}^2d_\Gamma^{\;2})(x)\geq0$ by Proposition~\ref{pcc2.1}.
In fact the lower bound is attained if $K$ has a proper face with dimension $d-1$.
\hfill$\Box$

\bigskip

At this point we are prepared to establish the weighted Hardy inequalities (\ref{ecc1.1}).

\smallskip

\noindent{\bf Proof of Theorem~\ref{tcc1.1}}$\;$
Let $\chi_p=c_\Omega\,d_\Gamma^{\;-p}\,(\nabla d_\Gamma^{\;2})$.
Further let $c'_\Omega=c'\circ d_\Gamma$.
Then 
\begin{eqnarray*}
\divv\chi_p&=&2\,(c'_\Omega \,d_\Gamma/c_\Omega-p)\,c_\Omega\,d_\Gamma^{\;-p}\,|\nabla d_\Gamma|^2
+c_\Omega\,d_\Gamma^{\;-p}\,(\nabla^2d_\Gamma^{\;2})\\[5pt]
&\geq& 2\,b_p\,c_\Omega\,d_\Gamma^{\;-p}
\end{eqnarray*}
with $b_p=(d-d_H+\delta\wedge\delta'-p)$ where we have used $|\nabla d_\Gamma|^2=1$,
the estimate $\nabla^2d_\Gamma^{\;2}\geq2\,(d-d_H)$ of Corollary~\ref{ccc2.2} and the observation that 
$c'_\Omega \,d_\Gamma/c_\Omega\geq \delta\wedge\delta'$.
(The last estimate follows since $s\,c'(s)/c(s)=(\delta+s\,\delta')/(1+s)$.)

Next for $\varepsilon>0$ set $\varphi_\varepsilon=(\varphi^2+\varepsilon^2)^{1/2}-\varepsilon$.
Then $\varphi_\varepsilon\geq0$ is a regularized approximation to $|\varphi|$ with the same support as $\varphi$.
But $\varphi^2+\varepsilon^2=(\varphi_\varepsilon+\varepsilon)^2\geq \varphi_\varepsilon^2+\varepsilon^2$
so $\varphi_\varepsilon\leq |\varphi|$.
In addition $\nabla\varphi_\varepsilon=(\varphi/(\varphi^2+\varepsilon^2)^{1/2})\,\nabla\varphi$.
Now assume $p\in\langle1,\infty\rangle$ and 
$b_p>0$.
Then
\begin{eqnarray*}
0<2\,b_p\int_\Omega c_\Omega\,d_\Gamma^{\;-p}\,\varphi_\varepsilon^{\,p}&\leq&\int_\Omega(\divv\chi_p)\,\varphi_\varepsilon^{\,p}\\[5pt]
&=&-p\int_\Omega c_\Omega\,d_\Gamma^{\;-p}\,(\nabla d_\Gamma^{\;2}).(\nabla\varphi_\varepsilon)\,\varphi_\varepsilon^{\,p-1}\\[5pt]
&=&-2\,p\int_\Omega c_\Omega\,d_\Gamma^{\;-p}\,(\nabla d_\Gamma).(\nabla\varphi_\varepsilon)\,\varphi_\varepsilon^{\,p-1}\\[5pt]
&\leq&2\,p\,\Big(\int_\Omega(c_\Omega\,d_\Gamma^{\;-p+1})^p|(\nabla d_\Gamma).(\nabla\varphi)|^p\,\psi^p\Big)^{1/p}.
\Big(\int_\Omega\varphi_\varepsilon^{\,p}\,\psi^{-q}\Big)^{1/q}
\end{eqnarray*}
where $q$ is the conjugate of $p$ and $\psi$ is a strictly positive function.
The last step uses the H\"older inequality.
Choosing  $\psi=c_\Omega^{-1/q}\,d_\Gamma^{\;p-1}$ one finds
\begin{eqnarray*}
0<b_p\int_\Omega c_\Omega\,d_\Gamma^{\;-p}\,\varphi_\varepsilon^{\,p}&\leq& p\,\Big(\int_\Omega c_\Omega\,|(\nabla d_\Gamma).(\nabla\varphi)|^p\Big)^{1/p}.
\Big(\int_\Omega c_\Omega\,d_\Gamma^{\;-p}\varphi_\varepsilon^{\,p}\Big)^{1/q}\;.
\end{eqnarray*}
Dividing by the last factor and raising the inequality to the power $p$ one obtains 
\[
\int_\Omega c_\Omega\,|\nabla\varphi|^p\geq \int_\Omega c_\Omega\,|(\nabla d_\Gamma).(\nabla\varphi)|^p\geq a_p\int_\Omega c_\Omega\,d_\Gamma^{\;-p}\,\varphi_\varepsilon^{\,p}
\]
for all $\varphi\in C_c^1(\Omega)$.
Then  the Hardy inequality of the proposition follows in the limit $\varepsilon\to0$ by dominated convergence.

The proof for $p=1$ is similar but simpler.
The H\"older inequality is not necessary.
\hfill$\Box$

\bigskip

The existence of a weighted Hardy inequality of the form (\ref{ecc1.1}) in the   situation, $\delta=\delta'$, and with $d_H<d-1$, follows from Theorem~4.2 of
\cite{LV}.
This paper also indicates a number of interesting directions to extend the current results.

\begin{remarkn}\label{rcc1} 
The foregoing proof only uses some general features of the weight function $c$.
The estimates (\ref{ecc1.1}) follow for any strictly positive differentiable $c$ on $\langle0,\infty\rangle$ with 
$c'(s)s/c(s)\geq \delta\wedge \delta'$.
If one makes the replacement $c(s)\to c(s)=s^\delta(a+bs)^{\delta'-\delta}$ with $a, b>0$ then 
$c'(s)s/c(s)=(a\,\delta+b\,\delta's)/(a+b\,s)\geq \delta\wedge \delta'$ and the theorem remains valid. 
Moreover, the constant $a_p$ in the Hardy inequality (\ref{ecc1.1})
is unchanged but now $c(s)/s^\delta\to a^{\delta'-\delta}$ as $s\to0$ and $c(s)/s^{\delta'}\to b^{\delta'-\delta}$ as $s\to\infty$. 
\end{remarkn}

\begin{remarkn}\label{rcc2} The condition  $d-d_H+\delta\wedge\delta'-p>0$ in  Theorem~\ref{tcc1.1} 
restricts the result to sets whose boundary have small codimension. 
For example if $\delta=0=\delta'$ it requires $d-d_H>p\geq1$.
In particular it does not apply if $d_H=d-1$.
If, however,  $d_H$  is small it is useful and allows one to  deduce Rellich inequalities on $L_2(\Omega)$ by the arguments
of \cite{Rob12} for all $\delta, \delta'\geq0$ (see Section~\ref{S3}).
\end{remarkn}

The foregoing arguments may also be used to obtain Hardy inequalities on convex subsets $\Omega$
but the conclusions are somewhat weaker.
The problem is that  points in $\Omega$ can have multiple near points.
This causes complications since $|(\nabla d_\Gamma)(x)|=1$ if and only if $x$ has a unique near point (see \cite{BEL}, Section~2.2).
The set of points in $\Omega$ which have more than one near point is defined as the skeleton $S(\Omega)$ of the set.
Then $ |(\nabla d_\Gamma)(x)|=1$ on $G(\Omega)=\Ri^d\backslash\overline{S(\Omega)}$.

The following result is in the spirit of Theorem 3.4.3 of \cite{BEL}.

\begin{prop}\label{p2}
Assume $\Omega$ is convex.
Again let $c_\Omega=c\circ d_\Gamma$ with $c(s)=s^\delta(1+s)^{\delta-\delta'}$ where $\delta,\delta'\geq0$.
If $p-1-\delta\vee\delta'>0$ then
\[
\int_\Omega c_\Omega\,|\nabla\varphi|^p\geq \int_\Omega c_\Omega\,|(\nabla d_\Gamma).(\nabla\varphi)|^p\geq a_p\int_\Omega c_\Omega\,d_\Gamma^{\;-p}\,|\varphi|^p
\]
for all $\varphi\in C_c^1(G(\Omega))$ with $a_p=((p-1-\delta\vee\delta')/p)^p$.
\end{prop}
\proof\
If $\Omega$ is  convex then $d_\Gamma$ is concave (see, for example, \cite{BEL}, Theorem~2.3.2).
This is sufficient to deduce that  $\Delta d_\Gamma$  is a positive measure
(see \cite{EvG}, Chapter~6).
Therefore
\[
\int_\Omega\,(\nabla\psi).(\nabla d_\Gamma)=\int_\Omega\,d\mu_\Omega\,\psi\geq0
\]
for all  positive $\psi\in C_c^1(\Omega)$ with $\mu_\Omega$  a positive Radon measure.
Again introduce the regularizations of $|\varphi|$ by
 $\varphi_\varepsilon=(\varphi^2+\varepsilon^2)^{1/2}-\varepsilon$ with   $\varepsilon>0$.
It then follows that 
\[
\int_\Omega\,(\nabla (c_\Omega\,d_\Gamma^{\;-p+1}\,\varphi_\varepsilon^{\,p})).(\nabla d_\Gamma)\geq0
\;.
\]
Therefore
\[
\int_\Omega\,(\nabla (c_\Omega\,d_\Gamma^{\;-p+1})).(\nabla d_\Gamma)\,\varphi_\varepsilon^{\,p}
+p\int_\Omega c_\Omega\,d_\Gamma^{\;-p+1}\,(\nabla d_\Gamma).(\nabla\varphi_\varepsilon)\,\varphi_\varepsilon^{\,p-1}
\geq0
\;.
\]
Next it follows that 
\begin{eqnarray*}
-(\nabla (c_\Omega\,d_\Gamma^{\;-p+1})).(\nabla d_\Gamma)&=&(p-1-c_\Omega'\,d_\Gamma/c_\Omega)\,c_\Omega\,d_\Gamma^{\;-p}
|\nabla d_\Gamma|^2\\[5pt]
&\geq&(p-1-(\delta\vee\delta'))\,c_\Omega\,d_\Gamma^{\;-p}|\nabla d_\Gamma|^2
\end{eqnarray*}
since $c_\Omega'\,d_\Gamma\leq (\delta\vee\delta')\,c_\Omega$.

Next if $\varphi\in C_c^1(G(\Omega))$ then $\supp\varphi_\varepsilon\subset G(\Omega)$ and consequently $|\nabla d_\Gamma|=1$
on the support of $\varphi_\varepsilon$.
Therefore, by combining the foregoing estimates , 
one obtains
\[
0<b_p\int_\Omega c_\Omega\,d_\Gamma^{\;-p}\,\varphi_\varepsilon^{\,p}\leq p\int_\Omega c_\Omega\,d_\Gamma^{\;-p+1}\,|(\nabla d_\Gamma).(\nabla\varphi)|\,\varphi_\varepsilon^{\,p-1}
\]
whenever $b_p=(p-1-(\delta\vee\delta'))>0$. Here we have again used $\nabla\varphi_\varepsilon=(\varphi/(\varphi^2+\varepsilon^2)^{1/2})\,\nabla\varphi$.
Therefore the H\"older inequality gives
\[
b_p\int_\Omega c_\Omega\,d_\Gamma^{\;-p}\,\varphi_\varepsilon^{\,p}
\leq p\,\Big(\int_\Omega(c_\Omega\,d_\Gamma^{\;-p+1})^p|(\nabla d_\Gamma).(\nabla\varphi)|^p\,\psi^p\Big)^{1/p}.
\Big(\int_\Omega\,\varphi_\varepsilon^{\,p}\,\psi^{-q}\Big)^{1/q}
\]
for all $\psi$ positive.
One can then proceed as previously and choose $\psi=c_\Omega^{-1/q}\,d_\Gamma^{\;p-1}$  to find
\begin{eqnarray*}
b_p\int_\Omega c_\Omega\,d_\Gamma^{\;-p}\,\varphi_\varepsilon^{\,p}&\leq &
p\,\Big(\int_\Omega c_\Omega\,|(\nabla d_\Gamma).(\nabla\varphi)|^p\Big)^{1/p}.
\Big(\int_\Omega c_\Omega\,d_\Gamma^{\;-p}\,\varphi_\varepsilon^{\,p}\Big)^{1/q}\\[5pt]
\end{eqnarray*}
Then since  $b_p>0$ one can divide throughout by the last factor, raise the inequality to the $p$-th power 
and take the limit $\varepsilon\to0$ to obtain the 
second inequality in the proposition.
The first one then follows since $|\nabla d_\Gamma|=1$ on $G(\Omega)$.
\hfill$\Box$
\bigskip

For further results on the weighted and unweighted Hardy inequality on convex sets we refer to \cite{MMP}, \cite{Avk1} and references therein.

\section{Rellich inequalities}\label{S3}

In this section we establish the  Rellich inequalities  (\ref{ecc1.2}) of Theorem~\ref{tcc1.2}.
Our proof is based on an extension of  Theorem~4 in the paper of Davies and Hinz \cite{DaH} (see Theorem~6.3.3 in  \cite{BEL})  from the Laplacian  to the weighted operator $H$.

\begin{prop}\label{pcc3.1}
Let $\Omega$ be a general domain in $\Ri^d$ and fix $p\in\langle1,\infty\rangle$.
Define the closeable operator $H=-\sum^d_{k=1}\partial_k \,c_\Omega\,\partial_k$ on  $D(H)=C_c^\infty(\Omega)$.
If there is a $\chi$ in the domain of the $L_p$-closure $\overline H$ of $H$ such that 
$\chi>0$, $\overline{H}\chi>0$ and $\overline{H}\chi^{1+\gamma}\geq0$ for some  $\gamma>0$ then
\begin{equation}
\int_\Omega|\overline{H}\chi|\,|\varphi|^p\leq p^{2p}(p+\gamma\,(p-1))^{-p} \int_\Omega \chi^p\,|\overline{H}\chi|^{-p+1}\,
|H\varphi|^p
\label{ecc3.1}
\end{equation}
for all $\varphi\in C_c^\infty(\Omega)$.
\end{prop}

This  proposition differs superficially from that of Davies--Hinz  since we define the Laplacian as $\Delta=-\nabla^2$ instead of $\nabla^2$.
Similarly we have  introduced a minus sign in the definition of $H$.
Moreover, the parameter $\delta$ in \cite{DaH} is replaced by $1+\gamma$ and this changes slightly  the form of the constant in (\ref{ecc3.1})

The proof of Proposition~\ref{pcc3.1} closely follows the arguments of \cite{DaH}.
The introduction of the coefficient $c_\Omega$ makes no essential change.
In fact since the estimates are on $C_c^\infty(\Omega)$  it suffices that $c_\Omega$ is the operator
of multiplication by a locally $C_1$-function.
The Davies--Hinz result also extends to more general divergence-form operators but this is not relevant in the current context.
It suffices that it applies to the weight functions used in Theorem~\ref{tcc1.2}.
Since the proof of Theorem~4 in \cite{DaH} is relatively long and since its adaptation to the weighted operators
does not introduce any significant changes we omit further discussion of the proof of Proposition~\ref{pcc3.1}.
We do, however, give the details of its application to the proof of Theorem~\ref{tcc1.2}.

\medskip

\noindent{\bf Proof of Theorem~\ref{tcc1.2}}$\;$
Define $\chi$ on  the open right half line by $\chi(s)=s^{-\alpha}(1+s)^{-\alpha'+\alpha}$ with $\alpha, \alpha'\geq0$.
Then set $\chi_\Omega=\chi\circ d_\Gamma$ and adopt the notation $\chi'_\Omega=\chi'\circ d_\Gamma$ etc.
Our aim is to derive conditions on $\alpha$ and $\alpha'$ such that $H\chi_\Omega>0$ with $H$ (the closure of)  the weighted operator of
Theorem~\ref{tcc1.2}.
In fact one can obtain quite precise lower bounds on $H\chi_\Omega$.

\begin{lemma}\label{lcc3.1} Let
$b_\alpha=\left(d-d_H+(\delta\wedge\delta')\right)(\alpha\wedge\alpha')
-(\alpha\vee\alpha')\,(\alpha\vee\alpha'+2)$.

\smallskip

It follows that  $H\chi_\Omega\geq b_\alpha\,d_\Gamma^{\;-2} \,c_\Omega\,\chi_\Omega$.  
Hence if $b_\alpha>0$ then $H\chi_\Omega>0$.
\end{lemma}
\proof\
First one has 
$\chi'(s)
=-s^{-1}\,\chi(s)\,(\alpha +\alpha' s)(1+s)^{-1}$.
Therefore
\[
-s^{-1}\,\chi(s)\,(\alpha\vee\alpha')\leq \chi'(s)\leq -s^{-1}\,\chi(s)\,(\alpha\wedge\alpha')
\;.
\]
In addition
\begin{eqnarray*}
\chi''(s)
&=&s^{-2}\,\chi(s)\,(1+s)^{-2}\,\Big(\alpha\,(\alpha+1)+2\,\alpha\,(\alpha'+1)\,s+\alpha'\,(\alpha'+1)\,s^2\Big)\\[5pt]
&\leq&s^{-2}\,\chi(s)\,(\alpha\vee\alpha')\,(\alpha\vee\alpha'+1)
\;.
\end{eqnarray*}
Secondly,  one calculates that 
\begin{eqnarray}
H\chi_\Omega&=&-d_\Gamma^{\;-1}c_\Omega\,\chi'_\Omega \left(\nabla^2d_\Gamma^{\;2}\right)/2-\left(c'_\Omega\,\chi'_\Omega
-d_\Gamma^{\;-1}\,c_\Omega\,\chi'_\Omega+c_\Omega\,\chi''_\Omega\right)|\nabla d_\Gamma|^2
\;.\label{ecc3.10}
\end{eqnarray}
But $|\nabla d_\Gamma|=1$ by the discussion at the beginning of Section~\ref{S2} and $(\nabla^2d_\Gamma^{\;2})/2\geq d-d_H$ by 
Corollary~\ref{ccc2.2}.
Then we use the bounds on $\chi'$ and $\chi''$ together with  the  lower bound  $c'(s)\geq (\delta\wedge\delta')\,s^{-1}c(s)$ to estimate  the four terms on the right hand side
of (\ref{ecc3.10}).
The first two  terms give positive contributions but  the other terms are negative.
One finds
\begin{eqnarray*}
H\chi_\Omega&\geq&
\left( (d-d_H)+(\delta\wedge\delta')\right)\,(\alpha\wedge\alpha') \left(d_\Gamma^{\;-2}\,c_\Omega\,\chi_\Omega\right)\\[4pt]
&&\hspace{4cm}{}-\left((\alpha\vee\alpha')+ (\alpha\vee\alpha')\,(\alpha\vee\alpha'+1)\right)\,\left(d_\Gamma^{\;-2}\,c_\Omega\,\chi_\Omega\right)\\[5pt]
&=&b_\alpha\, d_\Gamma^{\;-2}\,c_\Omega\,\chi_\Omega \;.
\end{eqnarray*}
 Clearly $H\chi_\Omega>0$ if the $\delta, \alpha,$ etc. are such that $b_\alpha>0$.
 \hfill$\Box$
 \bigskip
 
Now  assuming that $\alpha$ and $\alpha'$ are chosen to ensure that $b_\alpha>0$ one can bound the product
$\chi_\Omega^{\;p}\,|H\chi_\Omega|^{-p+1}$ occurring on the right hand side of (\ref{ecc3.1}).
Explicitly one obtains
\[
\chi_\Omega^{\;p}\,|H\chi_\Omega|^{-p+1}\leq b_\alpha^{-p+1}\,d_\Gamma^{\;-\sigma}(1+d_\Gamma)^{-\tau}
\]
with
$\sigma=\alpha-(2-\delta)(p-1)$ and 
$\tau=(\alpha'-\alpha)+(\delta'-\delta)(p-1)$.
Hence if one chooses $\alpha=\alpha_p=(2-\delta)(p-1)$ and $\alpha'=\alpha'_p=(2-\delta')(p-1)$
one obtains the uniform bound
\begin{equation}
\chi_\Omega^{\;p}\,|H\chi_\Omega|^{-p+1}\leq b_{\alpha_p}^{-p+1}
\label{ecc3.2}
\end{equation}
as long as 
\[
b_{\alpha_p}=\left(d-d_H+(\delta\wedge\delta')\right)(\alpha_p\wedge\alpha'_p)
-(\alpha_p\vee\alpha_p')\,(\alpha_p\vee\alpha_p'+2)>0
\;.
\]
But this is a  condition on $p, \delta$ and $ \delta'$.

\begin{lemma}\label{lcc3.2}
If $(d-d_H+p\,(\delta\wedge\delta')-2p)\geq 2p\,|\delta-\delta'|\,(2-\delta\vee\delta')^{-1}$ then $b_{\alpha_p}>0$.
\end{lemma}
\proof\
Substituting the values of $\alpha_p$ and $\alpha'_p$ in the definition of $b_\alpha$ one finds
\begin{eqnarray*}
b_{\alpha_p}
&=&\left(d-d_H+(\delta\wedge\delta')-(\alpha_p\vee\alpha_p'+2)\right)\,(\alpha_p\wedge\alpha_p')-|\alpha_p-\alpha_p'|\,((\alpha_p\vee\alpha_p')+2)\\[5pt]
&\geq&(p-1)\left((d-d_H+p\,(\delta\wedge\delta')-2\,p)(2-\delta\vee\delta')-2\,p\,|\delta-\delta'|\right)
\;.
\end{eqnarray*}
Since $p>1$ the statement of the lemma follows immediately.
\hfill$\Box$

\bigskip

Note that the condition of the lemma is  the  condition  posited in Theorem~\ref{tcc1.2} for validity of the Rellich inequality.

The next lemma  provides the last estimates necessary for the application of  Proposition~\ref{pcc3.1} to derive the Rellich inequality.
\begin{lemma}\label{lcc3.3}
Let $\tilde\chi_\Omega=d_\Gamma^{\;-\alpha_p}(1+d_\Gamma)^{-\alpha'_p+\alpha_p}$. Assume $b_{\alpha_p}>0$.
Then
\[
\tilde\chi_\Omega^{\;p}\,|H\tilde\chi_\Omega|^{-p+1}\leq b_{\alpha_p}^{\;-p+1}
\;\;\;\;\;\;\;and \;\;\;\;\;\;\; H\tilde\chi_\Omega\geq b_{\alpha_p}\,(c_\Omega\,d_\Gamma^{\;-1})^p
\;.
\]
Moreover, $ H\tilde\chi_\Omega^{\;1+\gamma}\geq0$ for all $\gamma\in[0,\gamma_p\,]$ where 
$\gamma_p=b_{\alpha_p}/(\alpha_p\vee\alpha'_p)^2$.
\end{lemma}
\proof\
The first estimate follows from Lemma~\ref{lcc3.1} and the choice of $\alpha_p$ and $\alpha'_p$ as discussed above.
The second estimate follows from another application of Lemma~\ref{lcc3.1} by noting that
\begin{eqnarray*}
H\tilde\chi_\Omega\geq b_{\alpha_p}\,d_\Gamma^{\;-2}\,c_\Omega\,\tilde\chi_\Omega
&=&b_{\alpha_p}\,d_\Gamma^{\;-2}\,d_\Gamma^{\;\delta}(1+d_\Gamma)^{\delta'-\delta}\,d_\Gamma^{\;-\alpha_p}(1+d_\Gamma)^{-(\alpha'_p-\alpha_p)}\nonumber\\[5pt]
&=&b_{\alpha_p}\,d_\Gamma^{\;-2p}\,d_\Gamma^{\;\delta p}(1+d_\Gamma)^{(\delta'-\delta)p}=b_{\alpha_p}\,(c_\Omega\,d_\Gamma^{\;-1})^p\end{eqnarray*}
where the second equality results from substituting the specific values of $\alpha_p$ and $\alpha_p'$.

The last statement of the lemma
follows by first noting that 
\[
\tilde\chi_\Omega^{\;1+\gamma}=d_\Gamma^{\;(1+\gamma)\alpha_p}(1+d_\Gamma)^{(1+\gamma)(-\alpha'_p+\alpha_p)}
\;.
\]
Therefore $ H\tilde\chi_\Omega^{\;1+\gamma}\geq0$ if $b_{(1+\gamma)\alpha_p}\geq0$ by a third application of Lemma~\ref{lcc3.1}.
But
\[
 b_{(1+\gamma)\alpha_p}=(1+\gamma)\,(b_{\alpha_p}-\gamma\,(\alpha_p\vee\alpha_p')^2)
\]
by the definition of $b_{\alpha}$.
Therefore $ b_{(1+\gamma)\alpha_p}\geq0$ whenever $0\leq\gamma\leq \gamma_p$.
\hfill$\Box$

\bigskip

At this point we have verified the  conditions necessary for the application of Proposition~\ref{pcc3.1} to $H$ and $\tilde\chi_\Omega$
to obtain the Rellich inequalities of Theorem~\ref{tcc1.2}.
We  now evaluate (\ref{ecc3.1}) with the foregoing estimates.
First we observe that $b_{\alpha_p}>0$ by Lemma~\ref{lcc3.2} and  the assumption of the theorem.
Then it follows from the estimates of Lemma~\ref{lcc3.3} that 
\begin{eqnarray*}
b_{\alpha_p}\int_\Omega |c_\Omega\,d_\Gamma^{\,-2}\varphi|^p&\leq&\int_\Omega |H\chi_\Omega|\,|\varphi|^p\\[5pt]
&\leq& p^{2p}(p+\gamma_p\,(p-1))^{-p} \int_\Omega \chi_\Omega^{\,p}\,|H\chi|^{-p+1}\,|H\varphi|^p\\[5pt]
&\leq& p^{2p}(p+\gamma_p\,(p-1))^{-p}\,b_{\alpha_p}^{-p+1}\int_\Omega |H\varphi|^p
\;.
\end{eqnarray*}
Thus  by rearrangement one obtains the Rellich inequality (\ref{ecc1.2}) with 
\[
c_p=(p+\gamma_p\,(p-1))\,b_{\alpha_p}\,p^{-2}
\;.
\]
It follows from $b_{\alpha_p}, \gamma_p>0$  that $c_p>0$.
We next  argue that $c_p\leq C_p$.

First one has
\begin{eqnarray*}
b_\alpha=\left(d-d_H+(\delta\wedge\delta')-(\alpha\vee\alpha'+2)\right) (\alpha\wedge\alpha')-a_\alpha
\end{eqnarray*}
with
\begin{eqnarray*}
a_\alpha&=&(\alpha\vee\alpha')\,(\alpha\vee\alpha'+2)-(\alpha\wedge\alpha')\,(\alpha\vee\alpha'+2)\\[5pt]
&=&|\alpha-\alpha'|\,((\alpha\vee\alpha')+2)\geq0
\;.
\end{eqnarray*}
Now set 
\[
\tilde b_\alpha=(d-d_H+(\delta\wedge\delta')-(\alpha\vee\alpha'+2))
\]
Then
\[
b_\alpha=(\alpha\wedge\alpha')\,\tilde b_\alpha-a_\alpha\leq (\alpha\wedge\alpha')\,\tilde b_\alpha
\;.
\]
Hence $b_{\alpha_p}\leq (\alpha_p\wedge\alpha'_p)\,\tilde b_{\alpha_p}$ with equality if and only if $\alpha_p=\alpha'_p$ or, equivalently,
$\delta=\delta'$.
Moreover, $\gamma_p=b_{\alpha_p}/(\alpha_p\vee\alpha'_p)^2\leq \tilde\gamma_p$,  where $\tilde\gamma_p=\tilde b_{\alpha_p}/(\alpha_p\vee\alpha'_p)$, with equality if and only if $\delta=\delta'$.
Now 
\begin{eqnarray*}
c_p&\leq&(p+\tilde\gamma_p\,(p-1))\,(\alpha_p\wedge\alpha'_p)\,\tilde b_{\alpha_p}\,p^{-2}\\[5pt]
&\leq&((\alpha_p\vee\alpha'_p)\,p+\tilde b_{\alpha_p}\,(p-1))\,\tilde b_{\alpha_p}\,p^{-2}
\;.
\end{eqnarray*}
But 
\begin{eqnarray*}
\tilde b_{\alpha_p}&=&(d-d_H+(\delta\wedge\delta')-((2-\delta\wedge\delta')\,(p-1)+2)\\[5pt]
&=&(d-d_H+p\,(\delta\wedge\delta')-2p)
\end{eqnarray*}
and
\begin{eqnarray*}
(\alpha_p\vee\alpha'_p)\,p+{\tilde b_{\alpha_p}}(p-1)&=&(2-(\delta\wedge\delta'))\,p(p-1)+{\tilde b_{\alpha_p}}(p-1)\\[5pt]
&=&(p-1)\,(d-d_H)
\end{eqnarray*}
Combining these estimates one has $c_p\leq C_p$ where $C_p$ is defined in Theorem~\ref{tcc1.2}.
\hfill$\Box$

\bigskip

We have avoided calculating $c_p$ explicitly since the resulting expression is complicated and is not necessarily optimal.
It is, however, straightforward to identify it from the value of $b_{\alpha_p}$ given prior to Lemma~\ref{lcc3.3} and the definition of
$\gamma_p$.
Nevertheless $c_p$ does have some simple properties as a function of the degeneracy parameters $\delta$ and $\delta'$.

Set $c_p=c_p(\delta,\delta')$ to denote the dependence on $\delta$ and $\delta'$.
Then $c_p$ is a positive symmetric function and $\delta\in[0,2\rangle\mapsto c_p(\delta,\delta)$ is strictly increasing.
Moreover,  if $c_p(\delta_0,0)\geq 0$ then $\delta\in[0,\delta_0]\mapsto c_p(\delta,0) $ is strictly decreasing.
In particular
\[
c_p(0,0)\geq c_p(\delta,0)\geq c_p(\delta_0,0)
\]
for all $\delta\in[0,\delta_0]$.

These inequalities follow because
\[
c_p(\delta,\delta)=(p-1)(d-d_H)(d-d_H+p\delta-2p)/p^2
\]
and 
\[
c_p(\delta,0)=c_p(0,\delta)=(p-1)(d-d_H)((d-d_H)(1-\delta/2)-2p)(1-\delta/2)/p^2
\]
which are special cases of the general formula for $c_p$.

\section{Optimal constants}\label{S4}

In this section we consider the problem of deriving optimal constants in the Hardy and Rellich inequalities of Theorems~\ref{tcc1.1}
and \ref{tcc1.2}.
First we discuss whether  the constant $a_p^{\,p}$ in Theorem~\ref{tcc1.1} is  the largest possible for the  Hardy inequality.
The maximal positive constant $\mu_p(\Omega)$ for which (\ref{ecc1.1}) is valid is given by
\begin{equation}
\mu_p(\Omega)=
\inf\Big\{\int_\Omega c_\Omega\,|\nabla\varphi|^p\Big/\!\int_\Omega c_\Omega\,d_\Gamma^{\;-p}\,|\varphi|^p:\; \varphi\in C_c^1(\Omega)\Big\}
\;.\label{ecc4.1}
\end{equation}
Clearly $\mu_p(\Omega)\geq a_p^{\,p}$ by Theorem~\ref{tcc1.1}. 
Therefore optimality follows if the infimum in (\ref{ecc4.1}) is less than or equal to  $a_p^{\,p}$.
Since $c_\Omega$ has a different asymptotic behaviour at the boundary $\Gamma$ to that at infinity this variational problem has two distinct elements, a local and a global.
For orientation we begin with a brief discussion of the classical case $K=\{0\}$ (see, for example, \cite{BEL} Section~1.2).

If $\Omega=\Ri^d\backslash\{0\}$ then the constant $a_p$ of Theorem~\ref{tcc1.1} is given by $a_p=(d+\delta\wedge\delta'-p)/p$.
Therefore   if  $a_p(\sigma)=((d+\sigma -p)/p)^p$ for all $\sigma\geq0$ then $a_p^{\,p}= a_p(\delta\wedge\delta')=a_p(\delta)\wedge a_p(\delta')$.
Thus to prove that $a_p^{\,p}=\mu_p(\Omega)$ it suffices to prove that $\mu_p(\Omega)\leq a_p(\delta)$ and $\mu_p(\Omega)\leq a_p(\delta')$.
This can be achieved by standard arguments  (see, for example, \cite{BEL}, Chapter~1).

The first upper bound follows by estimating the infimum in (\ref{ecc4.1})
 with a sequence of   functions $\varphi_\alpha=d_\Gamma^{\;-\alpha}\xi$, $\alpha>0$, where $\xi$ has support in  a small neighbourhood of the origin.
Since $d_\Gamma(x)=|x|$ it follows that  $c_\Omega\,|\nabla\varphi|^p$ and $c_\Omega\,d_\Gamma^{\;-p}\,|\varphi_\alpha|^p$ are integrable at the origin if $\alpha <(d+\delta-p)/p$.
In which case the leading term $d_\Gamma^{\;-\alpha}$ gives a bound proportional to $\alpha^p$ in the evaluation of  (\ref{ecc4.1}).
Then by a suitable choice of  localization functions $\xi$ and a  limiting argument one concludes that 
$\mu_p(\Omega)\leq ((d+\delta-p)/p)^p=a_p(\delta)$.
Here the property $\lim_{s\to0}c(s)\,s^{-\delta}=1$ is important.
The estimate at infinity is similar. 
One now chooses $\varphi_\alpha$ with support in the complement of  a large ball centred at the origin and again proportional to $d_\Gamma^{\,-\alpha}$.
Then, however,  $c_\Omega\,|\nabla\varphi|^p$ and $c_\Omega\,d_\Gamma^{\;-p}\,|\varphi_\alpha|^p$ are integrable at infinity if $\alpha>(d+\delta'-p)/p$.
Again the leading term gives  a bound proportional to $\alpha^p$.
Then another approximation and limiting argument gives the 
 upper bound
$\mu_p(\Omega)\leq ((d+\delta'-p)/p)^p=a_p(\delta')$.
Here the property $\lim_{s\to\infty}c(s)\,s^{-\delta'}=1$ is crucial.
Thus one arrives at the following conclusion.

\begin{prop}\label{pcc4.1}
If $K=\{0\}$ then the optimal constant $\mu_p(\Omega)$   in the Hardy inequality $(\ref{ecc1.1})$ is given by $\mu_p(\Omega)=a_p^{\,p}=((d+\delta\wedge\delta'-p)/p)^p$.
\end{prop}

In the more general situation that $\dim(K)\geq 1$  the foregoing approach is complicated by the geometry. 
Nevertheless one can obtain bounds by  a similar two step process of local estimates and  estimates at infinity.
The local estimates  are obtained by the methods of Barbatis, Filippas and Tertakis \cite{BFT}
which are also developed in Section~5 of Ward's thesis \cite{War}.
The following theorem covers the cases with $\dim(K)<d$.

\begin{thm}\label{tcc4.1} Adopt the assumptions of Theorem~$\ref{tcc1.1}$.
Further assume that  $ \dim(K)\in\{1,\dots, d-1\}$.

 Then the optimal constant $\mu_p(\Omega)$  in $(\ref{ecc1.1})$ satisfies $\mu_p(\Omega)\leq ((d-d_H+\delta-p)/p)^p$.
In particular, if $\delta\leq\delta'$ then $\mu_p(\Omega)= ((d-d_H+\delta-p)/p)^p$.
\end{thm}
\proof\
The  proposition   follows by the proof of Theorem~5.2.1 in \cite{War} but with some modification to take into account the weighting factor $c_\Omega$.
We outline a variation of Ward's argument which is subsequently extended to give a local bound on the optimal constant in the 
Rellich inequality (\ref{ecc1.2}).
First we give the proof for the special case $\delta=\delta'$ or, equivalently, $c(s)=s^\delta$.
Then, since the argument only involves functions with support in an arbitrarily small ball centred at a point of the boundary 
the result can be extended to the general  weighting factor $c_\Omega$.

The starting point of the  proof is a modification of Ward's Lemma~5.1.1.

\begin{lemma}\label{lcc4.1} Assume $c(s)=s^\delta$ with $\delta\geq0$.
Then
\begin{eqnarray}
\mu_p(\Omega)\leq(1-\lambda)^{-(p-1)} |(\beta+\delta-p)/p|^p
+\lambda^{-(p-1)}
\Big(\int_\Omega \,d_\Gamma^{\;-\beta+p} \,|(\nabla\varphi)|^p\Big/\!
\int_\Omega \,d_\Gamma^{\;-\beta} \,|\varphi|^p\Big)\label{ecc4.11}
\end{eqnarray}
for all $\varphi\in W^{1,p}_0(\Omega)$, $\beta\geq0$,  $p>1$ and $\lambda\in\langle0,1\rangle$.
\end{lemma}
\proof\ The proof follows that of Ward with $\varphi$  in (\ref{ecc4.1}) replaced by $\psi\varphi $
where $\psi=d_\Gamma^{\;-(\beta+\delta-p)/p}$.
Then  one uses the Leibniz rule, the 
$l_2$-triangle inequality
\[
|(\nabla\psi)\, \varphi+\psi\,(\nabla\varphi)|\leq |(\nabla\psi)|\, |\varphi|+|\psi|\,|(\nabla\varphi)|
\]
and the estimate
\begin{equation}
(s+t)^{p}\leq (1-\lambda)^{-(p-1)}s^p+\lambda^{-(p-1)}t^p
\label{ecc4.12}
\end{equation}
which is valid for all $s,t\geq0$, $\lambda\in\langle0,1\rangle$ and $p>1$.
(The latter inequality 
was used by  Secchi, Smets and Willem \cite{SSW} in their analysis of the Hardy inequality
on the complement of affine subsets.
It follows by minimization of the right hand side over $\lambda$.)
Now by combination of these observations one finds
\begin{eqnarray*}
\int_\Omega d_\Gamma^{\;\delta}\,|\nabla(\psi\varphi)|^p&\leq&
(1-\lambda)^{-(p-1)}\int_\Omega d_\Gamma^{\;\delta}\,|\nabla\psi|^p|\varphi|^p
+\lambda^{-(p-1)}\int_\Omega d_\Gamma^{\;\delta}\,|\psi|^p|\nabla\varphi|^p\\[5pt]
&=&(1-\lambda)^{-(p-1)}|(\beta+\delta-p)/p|^p
\int_\Omega d_\Gamma^{\;-\beta}\,|\varphi|^p+\lambda^{-(p-1)}\int_\Omega d_\Gamma^{\;-\beta+p}\,|\nabla\varphi|^p
\end{eqnarray*}
where we have used the explicit form of $\psi$.
Similarly
\[
\int_\Omega d_\Gamma^{\;\delta-p}\,|\psi\varphi|^p=\int_\Omega d_\Gamma^{\;-\beta}|\varphi|^p
\;.
\]
The statement of the lemma follows immediately.
\hfill$\Box$

\bigskip

The estimate for  $\mu_p(\Omega)$ given in Theorem~\ref{tcc4.1} now  follows by Ward's reasoning in  the proof of his Theorem~5.2.1.
The idea is to construct a sequence of $\varphi_n $ such that the numerator in the last term in (\ref{ecc4.11}) is bounded uniformly in $n$ if  $\beta=d-k$,
with $k=\dim(\partial K)=\dim(\Gamma)=d_H$,
but the denominator diverges as $n\to\infty$.
This is particularly easy in the current context since we are assuming $k=\dim(K)\leq d-1$.

First let $\Ri^d=\Ri^k\times \Ri^{d-k}$  where $\Ri^k$ is identified with the affine hull of $K$.
Therefore if  one sets  $x=(y,z)\in \Omega$ with $y\in \Ri^k$ and $z\in\Ri^{d-k}$ then 
$d_{\Omega}(x)=(d_K(y)^2+|z|^2)^{1/2}$ where $d_K(y)=\inf_{y'\in K}|y-y'|$.
Since $d_K(y)=0$ if $y\in K$ it follows that $d_\Gamma(y,z)=|z|$ if $y\in K$.
Secondly, define $\varphi\in C_c^\infty(\Omega)$ by setting $\varphi(y,z)=\eta(y)\chi(z)$ where $\eta\in C_c^\infty(K)$ and $\chi\in C_c^\infty(\Ri^{d-k}\backslash\{0\})$.
Further assume $\chi$ is a radial function.
Then with $\beta=d-k$ one has 
\begin{equation}
\int_\Omega\,d_\Gamma^{\,-\beta}|\varphi|^p=\int_Kdy\,|\eta(y)|^p\int_{\Ri^{d-k}}dz\,|z|^{-(d-k)}|\chi(z)|^p
=a_1\int_0^\infty dr\,r^{-1}|\chi(r)|^p
\label{ecc4.81}
\end{equation}
but 
\begin{eqnarray}
\int_\Omega\,d_\Gamma^{\,-\beta+p}|\nabla\varphi|^p&=&
\int_Kdy\int_{\Ri^{d-k}}dz\,|z|^{\,-(d-k-p)}(|(\nabla\eta)(y)|\,|\chi(z)|+|\eta(y)|\,|(\nabla\chi)(z)|)^p\nonumber\\[5pt]
&\leq&a_2\int_0^\infty dr\,r^{p-1}|\chi(r)|^p+a_3\int_0^\infty dr\,r^{p-1}|\chi'(r)|^p
\label{ecc4.82}
\end{eqnarray}
with $a_1, a_2, a_3>0$.
The last estimate  again uses (\ref{ecc4.12}).

Next consider the sequence of functions $\xi_n$ defined on $\langle0,\infty\rangle$ by $\xi_n(r)=0$ if $r\leq n^{-1}$,
$\xi_n(r)=\log rn/\log n$ if $n^{-1}\leq r\leq1$ and $\xi_n=1$ if $r\geq1$.
Then $0\leq \xi_n\leq1$ and the $\xi_n$ converge monotonically upward to the identity function.
Further let $\zeta$ be a $C^\infty$-function with $\zeta(r)=1$ if $r\leq1$, $\zeta(r)=0$ if $r\geq 2$
and $0\leq \zeta\leq1$.
Then set $\chi_n=\xi_n\zeta$.

It follows immediately that 
\[
\lim_{n\to\infty}\int_0^\infty dr\,r^{-1}|\chi_n(r)|^p=\infty
\;.
\]
Moreover,
\[
\int_0^\infty dr\,r^{p-1}|\chi_n(r)|^p\leq \int_0^2 dr\,r^{p-1}|\xi_n(r)|^p\leq \int_0^2 dr\,r^{p-1}\leq 2^{\,p}p^{-1}
\]
for all $ n>1$.
But   $\supp\chi_n\subseteq [0,2]$,  $\chi'_n=\xi'_n$ on $\langle0,1]$ and $ \chi'_n=\zeta'$ on $[1,2]$.
Therefore 
\[
 \int_0^\infty dr\,r^{p-1}|\chi_n'(r)|^p= \int_0^1 dr\,r^{p-1}|\xi_n'(r)|^p+ \int_1^2 dr\,r^{p-1}|\zeta'(r)|^p
 =(\log n)^{-(p-1)}+a
 \]
 where $a>0$ is the contribution of the second integral.
 Since $p>1$ the bound is uniform for all $n>1$.
 Hence if one sets $\varphi_n=\eta\,\chi_n$ one deduces from (\ref{ecc4.81}) and (\ref{ecc4.82}),
  with $\chi$ replaced by $\chi_n$, that 
  \[
 \limsup_{n\to\infty}\int_\Omega\,d_\Gamma^{\,-\beta+p}|\nabla\varphi_n|^p\Big/ \int_\Omega\,d_\Gamma^{\,-\beta}|\varphi_n|^p=0
 \;.
 \]
Therefore replacing  $\varphi$ with  $\varphi_n$ in (\ref{ecc4.11}) and setting $\beta=d-k$ one deduces that
\[
\mu_p(\Omega)\leq (1-\lambda)^{-(p-1)}\,((d-k+\delta-p)/p)^p
\]
for all $\lambda\in\langle0,1\rangle$.
Thus in the limit $\lambda\to0$ one has $\mu_p(\Omega)\leq a_p(\delta)$.

This completes the proof of the upper bound for $c(s)=s^\delta$,  i.e.\ for $\delta=\delta'$.
But it follows by construction that 
\[
\supp\varphi_n\subseteq \{(y,z): y\in \supp\eta, |z|\leq2\}
\;.
\]
The choice of the value $2$ is, however, arbitrary and by rescaling the $\xi_n$ it can be replaced  by any
$r>0$ without materially affecting the argument.
Then since $|z|^\delta(1+r)^{-|\delta-\delta'|}\leq c(z)\leq |z|^\delta(1+r)^{|\delta-\delta'|}$ for $|z|<r$  the case of general
$c_\Omega$ is reduced to the special case $\delta=\delta'$.

Finally if $\delta\leq \delta'$ it follows from Theorem~\ref{tcc1.1} that  $\mu_p(\Omega)\geq a_p(\delta)$.
Consequently one must have equality.
\hfill$\Box$

\bigskip

Next we investigate the derivation of the  bounds $\mu_p(\Omega)\leq ((d-d_H+\delta'-p)/p)^p$
in the setting of Theorem~\ref{tcc4.1}.
These bounds require additional information on the global properties of $K$.

The dimension $k$ of the convex set is defined as the dimension of the affine hull $A_K$ of $K$
and is essentially a local concept. 
It carries little information about the global character of the set.
For example, in $2$-dimensions $K$ could be a disc, an infinitely extended strip or a quadrant.
But viewed from afar these sets would appear to have dimension $0$, $1$ and $2$, respectively.
This aspect of the sets is captured by the `dimension at infinity' $k_\infty$ which is defined by
\[
k_\infty=\liminf_{r\to\infty}\left(\log|K\cap B_r|/\log r\right)
\]
where $B_r=\{y\in \Ri^k:|y|<r\}$ and $|S|$ indicates the $k$-dimensional Lebesgue measure of 
the set $S$.
The parameter $k_\infty$  of the convex set is integer valued  with  $0\leq k_\infty\leq k$.
In the two-dimensional examples it  takes the values $0$, $1$ and $2$  as expected.
The equality $k_\infty=k$ of the global and local dimensions will be the key property in deriving the upper bounds on $\mu_p(\Omega)$. 

\begin{lemma}\label{lcc4.0}
Assume   $k_\infty=k$.
Then 
\[
\inf_{\eta\in C_c^\infty(K)}\int_K|\nabla\eta|^p\Big/\!\!\int_K|\eta|^p=0
\;.
\]
\end{lemma}
\proof\
First let $\xi$ be a $C^\infty$-function with $0\leq\xi\leq 1$ such that 
$\supp\xi\subseteq K$ and $\xi(y)=1$ if $d_K(y)\geq1$, with $d_K$ the Euclidean distance to the boundary $\partial K$.
Secondly, let $\zeta_n$  be a sequence of $C^\infty$-functions
with $0\leq\zeta_n\leq 1$, $\zeta_n(y)$ if $y\in B_r$ and $\zeta_n=0$ if $y\in B^c_{r+1}$.
We may assume $\sup_n|\nabla\zeta_n|<\infty$.
Now set $\eta_n=\zeta_n\,\xi$.
Then $\eta_n\in C_c^\infty(K)$ and $\supp|\nabla\eta_n|$ has measure at most $b\,r^{k-1}$ for all $r\geq1$ with
$b>0$ independent of $r$.
But $\eta_n=1$ on a set of measure $c\,r^k$ with $c>0$.
Therefore 
\[
\int_K|\nabla\eta_n|^p\Big/\!\int_K|\eta_n|^p< a\,r^{-1}
\]
with $a>0$ independent of $r$.
The lemma follows immediately.
\hfill$\Box$

\bigskip

The following theorem establishes that $k_\infty=k$ is a sufficient condition for the expected global bounds but it is likely that it, or some variation of it, is also necessary.

\begin{thm}\label{tcc4.3}
Let $K$ be a closed convex subset of  $\Ri^d$ with $k=\dim(K)\in \{1,\ldots, d-1\}$ and with $k_\infty=k$.
Then the optimal constant $\mu_p(\Omega)$ in the Hardy inequality $(\ref{ecc1.1})$
on $\Omega=\Ri^d\backslash K$ is given by
$\mu_p(\Omega)=((d-k+\delta\wedge \delta'-p)/p)^p$.
\end{thm}
\proof\
First, $\mu_p(\Omega)\geq a_p^{\,p}$ with $a_p=(d-k+\delta\wedge \delta'-p)/p$ by  Theorem~\ref{tcc1.1}.
Therefore it suffices to establish a matching upper bound.
But   the local estimates of Theorem~\ref{tcc4.1} give the bound $\mu_p(\Omega)\leq ((d-k+\delta-p)/p)^p$.
Thus it remains to prove that $\mu_p(\Omega)\leq ((d-k+\delta'-p)/p)^p$.

Secondly, we again consider the decomposition $\Ri^d=\Ri^k\times \Ri^{d-k}$ with $K\subseteq \Ri^k
$ and $\Ri^k=A_K$.
Then since $d_\Gamma(y,z)=|z|$ if $y\in K$ the weighted Hardy inequality (\ref{ecc1.1}) on $L_p(\Omega)$  takes the form
\begin{equation}
\int_{K}dy\int_{\Ri^{d-k}} dz\,c(|z|)|(\nabla\varphi)(y,z)|^p\geq 
 a_p^{\,p}\int_{K}dy\int_{\Ri^{d-k}} dz\, c(|z|)|z|^{-p}|\varphi(y,z)|^p
\label{ecc4.2}
\end{equation}
for all $\varphi\in C_c^1(\Omega)$ with $\supp\varphi\subseteq K\times\Ri^{d-k}$ where
$c(s)=s^\delta(1+s)^{\delta'-\delta}$.
Therefore the  optimal constant  satisfies
\[
\mu_p(\Omega)\leq
\bigg(\int_{\Ri^k}dy\int_{\Ri^{d-k}} dz\,c(|z|)|(\nabla\varphi)(y,z)|^p\Big/
\!\int_{\Ri^k}dy\int_{\Ri^{d-k}} dz\,  c(|z|)|z|^{-p}|\varphi(y,z)|^p\bigg)
\;.
\]
Again let $\varphi$ be a product $\varphi(y,z)=\eta(y)\chi(z)$ with $\eta\in C_c^\infty(K)$ but  $\chi\in C_c^\infty(O_{\!R})$ where $O_{\!R}=\{z\in \Ri^{d-k}:|z|>R\}$.
Then
\[
\mu_p(\Omega)\leq{{ \int_{O_{\!R}} dz\,c(|z|)\int_K dy\Big(|(\nabla\chi)(z)|\,|\eta(y)|+|\chi(z)|\,|(\nabla\eta)(y)|\Big)^{p}}\over{
\Big(\int_{O_{\!R}}dz\, c(|z|)|z|^{-p}|\chi(z)|^p}\Big)\Big(\int_Kdy\,|\eta(y)|^p\Big)}
\;.
\]
We can again use (\ref{ecc4.12}) to estimate the right hand side.
One immediately obtains
\begin{eqnarray}
\mu_p(\Omega)&\leq&(1-\lambda)^{-(p-1)}\Bigg({{\int_{O_{\!R}} dz\,c(|z|)|(\nabla\chi)(z)|^p}\over{ \int_{O_{\!R}}dz\,c(|z|)|z|^{-p}|\chi(z)|^p}}\Bigg)
\nonumber\\[5pt]
&&\hspace{3.3cm}+\lambda^{-(p-1)}\Bigg({{\int_{O_{\!R}} dz\,c(|z|)|\chi(z)|^p}\over{ \int_{O_{\!R}} dz\,c(|z|)|z|^{-p}|\chi(z)|^p}}\Bigg)
\Bigg({{\int_K dy\,|(\nabla\eta)(y)|^p}\over{ \int_K dy\,|\eta(y)|^p}}\Bigg)
\label{ecc4.14}
\end{eqnarray}
for all $\lambda\in\langle0,1\rangle$.
Therefore
taking the infimum over $\eta\in C_c^\infty(K)$ 
followed by the infimum over $\lambda\in\langle0,1\rangle$ one deduces from Lemma~\ref{lcc4.0} that
\begin{equation}
\mu_p(\Omega)\leq 
{{\int_{O_{\!R}} dz\,c(|z|)|(\nabla\chi)(z)|^p}\over{ \int_{O_{\!R}}dz\,c(|z|)|z|^{-p}|\chi(z)|^p}}
\label{ecc4.13}
\end{equation}
for all $\chi\in C_c^\infty(O_{\!R})$ and all large $R$.

Finally  the infimum of the right hand side of (\ref{ecc4.13})  over $\chi$ followed by the limit $R\to\infty$  gives 
$\mu_p(\Omega)\leq ((d-k+\delta'-p)/p)^p$
by the global estimates for the Hardy inequality on $\Ri^{d-k}\backslash\{0\}$ sketched at the beginning of the
section.
The proof of the theorem now follows from this estimate combined with the  observations in the first paragraph of the proof.
\hfill$\Box$

\bigskip

Theorem~\ref{tcc4.3} applies to the special case that $K$ is an affine set since the assumption $k_\infty=k$ is automatically fulfilled.
The corresponding statement is an extension of a  result of \cite{SSW}.
Moreover, if $K$ is a general closed convex set and $A_K$ its  affine hull then  the theorem identifies the constant $a_p^{\,p}$ of Theorem~\ref{tcc1.1} as the optimal
constant $\mu_p(\Ri^d\backslash A_K)$ of the Hardy inequality (\ref{ecc1.1}) on $L_p(\Ri^d\backslash A_K)$.
Therefore one has the general conclusion that  $\mu_p(\Ri^d\backslash A_K)\leq  \mu_p(\Ri^d\backslash K)$
for convex sets with $\dim(K)=k\in\{1,\ldots,d-1\}$.
Moreover, $ \mu_p(\Ri^d\backslash A_K)=\mu_p(\Ri^d\backslash K)$ if $\delta\leq\delta'$ because the proof only requires a local estimate.

\medskip

Next we address the question of calculating the optimal constant in the Rellich inequality (\ref{ecc1.2}),
i.e.\ the value of 
\begin{equation}
\nu_p(\Omega)=
\inf\Big\{\int_\Omega \,|H\varphi|^p\Big/\!\int_\Omega c_\Omega^{\;p}\,d_\Gamma^{\;-2p}|\varphi|^p:\; \varphi\in C_c^2(\Omega)\Big\}
\;.\label{ecc4.4}
\end{equation}
Theorem~\ref{tcc1.2}  gives the lower bound  $\nu_p(\Omega)\geq c_p^{\,p}$ but this is  rather complicated and not likely to be
an efficient bound in general.
Therefore we consider the special case  $\delta=\delta'$ with weighting factor $d_\Gamma^{\,\delta}$.
Then Theorem~\ref{tcc1.2}  gives the simpler  bound $\nu_p(\Omega)\geq C_p^{\,p}$  with 
$C_p=(p-1)(d-d_H)(d-d_H+p\,\delta-2p)p^{-2}$.
Now we establish that $C_p^{\,p}$ is the optimal constant if $\delta=\delta'$ and  $\dim{K}\in \{0,1,\ldots,d-1\}$.
First we consider the degenerate case.

\begin{prop}\label{pcc4.2}
If $K=\{0\}$ and $\delta=\delta'\in[0,2\rangle$ then the optimal constant in the Rellich inequality $(\ref{ecc1.2})$ is given by
\[
\nu_p(\Omega)=C_p^{\,p}=\left((p-1)\,d\,(d+p\,\delta-2p)\,p^{-2}\right)^p
\]
for all $p>1$ for which $d+p\,\delta-2p>0$.
\end{prop}
\proof\
It follows from Theorem~\ref{tcc1.2}, with $\delta=\delta'$, that the lower bound $\nu_p(\Omega)\geq C_p^{\,p}$ is valid.
Therefore it suffices to establish a matching upper bound.
This is well known if $\delta=0$ but the proof is almost identical for $\delta\neq0$.
First, since $K=\{0\}$ one has $d_\Gamma(x)=|x|$.
Then as $\delta=\delta'$ one can deduce  an upper bound from (\ref{ecc4.4}) by a local estimate (see, for example, \cite{BEL}
Corollary~6.3.5 for the case of the Laplacian).
This is  achieved by the elementary procedure used to estimate the upper bound on the Rellich constant 
in the one-dimensional case.
One estimates with  radial functions $\varphi(x)=|x|^{-\alpha}\,\chi(|x|)$
where $\alpha>0$ and $\chi$ is a $C^2$-function with compact support near the origin.
The integrability of $|H\varphi|^p$ at the origin
imposes the restriction 
$d+p\,\delta-2p>0$.
Therefore one chooses $\alpha=(d+p\,\delta-2p+\varepsilon)/p$, with $\varepsilon>0$,
and estimates as in the one-dimensional case.
This leads to the upper bound $\nu_p(\Omega)\leq C_p^{\,p}$.
We omit the details.
\hfill$\Box$


\begin{remarkn}\label{rcc4.1}  If $K=\{0\}$ and $\delta\neq\delta'$ then one can establish the  upper bound  
$\nu_p(\Omega)\leq \left((p-1)\,d\,(d+p\,(\delta\wedge\delta')-2p)\,p^{-2}\right)^p$.
This follows by a local estimate which gives the bound $\left((p-1)\,d\,(d+p\,\delta-2p)\,p^{-2}\right)^p$ followed
by a similar estimate at `infinity' which gives the bound $\left((p-1)\,d\,(d+p\,\delta'-2p)\,p^{-2}\right)^p$.
Then one takes the minimum of the two bounds.
Unfortunately Theorem~\ref{tcc1.2} only gives a matching lower bound if $\delta=\delta'$.
If, for example, $\delta'=0$ then the upper bound is equal to $c_p(0,0)^p=\left((p-1)\,d\,(d-2p)\,p^{-2}\right)^p$
where we have used the notation introduced at the end of Section~\ref{S3}.
But Theorem~\ref{tcc1.2}  gives the lower bound $c_p(\delta,0)^p$ under the assumption that $c_p(\delta,0)>0$.
It follows, however, that  $c_p(\delta,0)<c_p(0,0)$ if $\delta>0$ by the discussion in Section~\ref{S3}.
\end{remarkn}

Now we establish a similar conclusion for  $\dim(K)\in\{1,\ldots,d-1\}$.
The following result corresponds to the Rellich analogue of Theorem~\ref{tcc4.1} and Theorem~\ref{tcc4.3}. 

\begin{thm}\label{tcc4.31}
Let $K$ be a closed convex subset of  $\Ri^d$ with $k=\dim(K)\in \{1,\ldots, d-1\}$ 
Then the optimal constant in the Rellich  inequality $(\ref{ecc1.2})$
satisfies the upper bound
\[
\nu_p(\Omega)\leq \left((p-1)(d-k)(d-k+p\,\delta-2p)\,p^{-2}\right)^p
\;.
\]
If, in addition, $k_\infty=k$ then
\[
\nu_p(\Omega)\leq \left((p-1)(d-k)(d-k+p\,(\delta\wedge\delta')-2p)\,p^{-2}\right)^p
\]
and for $\delta=\delta'$ one has equality.
\end{thm}
\proof\
The proof follows the earlier two step process of obtaining a local bound, dependent on $\delta$,
followed by a global bound, dependent on $\delta'$.
The local bound is independent of the assumption $k_\infty=k$.

\smallskip

\noindent{\bf Step 1}$\;$ The first statement of the theorem is established by a generalization of the 
local estimates used to prove Theorem~\ref{tcc4.1}.
Since all the estimates in this first step are local we again assume initially that $c(s)=s^\delta$.

Following  the earlier proof we choose coordinates   $x=(y,z)\in \Omega$ with $y\in \Ri^k$ and $z\in\Ri^{d-k}$
where  $\Ri^k$ is identified with the affine hull of $K$.
Then  $d_\Gamma(y,z)=|z|$ if $y\in K$.
Again we  define $\varphi\in C_c^\infty(\Omega)$ by setting $\varphi(y,z)=\eta(y)\chi(z)$ where $\eta\in C_c^\infty(K)$ and $\chi\in C_c^\infty(\Ri^{d-k}\backslash\{0\})$ is a radial function.
Next for $\alpha\geq0$ we set $\varphi_\alpha=d_\Gamma^{\;-\alpha}\varphi=\eta\,\chi_\alpha$ 
where $\chi_\alpha(z)=|z|^{-\alpha}\chi(z)$.
Thus $\varphi_\alpha=d_0^{\;-\alpha}\varphi$ where $d_0$ is the operator of multiplication by $|z|$.
Then 
\[
H\varphi_\alpha=(Hd_0^{\,-\alpha})\varphi+d_0^{\,-\alpha}(H\varphi)+2\,d_0^{\,\delta}(\nabla d_0^{\,-\alpha}).(\nabla\varphi)
\;.
\]
Therefore one calculates that 
\[
|H\varphi_\alpha|\leq \alpha(d-k+\delta-\alpha-2)d_0^{\,-\alpha-2+\delta}|\varphi|+R_\alpha
\]
if $d-k+\delta-2>\alpha$ where
\[
R_\alpha=d_0^{\,-\alpha}|H\varphi|+2\,\alpha\, d_0^{\,-\alpha-1+\delta}|\nabla\varphi|
\;.
\]
Hence it follows as in the proof of Lemma~\ref{lcc4.1} that
\begin{equation}
\nu_p(\Omega)\leq (1-\lambda)^{-(p-1)}( \alpha(d-k+\delta-\alpha-2))^p
+\lambda^{-(p-1)}\Big(\int_\Omega|R_\alpha|^p\Big/\int_\Omega d_0^{\,p(\delta-2)}|\varphi_\alpha|^p\Big)
\label{recc4.0}
\end{equation}
for all $\lambda\in\langle0,1\rangle$. 
Now we choose  $\alpha=(d-k+p\,\delta-2p)/p$ and assume $\alpha>0$.
Then  the constant appearing in the first term on the right is 
$\left((p-1)(d-k)(d-k+p\,\delta-2p)\,p^{-2}\right)^p$.
So it remains to prove that the second term, with the specific  choice of $\alpha$, can be made 
insignificant by a suitable choice of a sequence of  $\chi$.
First one has
\begin{equation}
\int_\Omega d_0^{\,p(\delta-2)}|\varphi_\alpha|^p=\int_Kdy\,|\eta(y)|^p\int_{\Ri^{d-k}}dz\,|z|^{-p(\alpha-\delta+2)}|\chi(z)|^p
=a_1\int^\infty_0dr\,r^{-1}|\chi(r)|^p
\label{recc4.1}
\end{equation}
with $a_1>0$.
Secondly,
\begin{eqnarray*}
|R_\alpha|^p&\leq& a\Big(d_0^{\,-p\alpha}|H\varphi|^p+d_0^{\,-p(\alpha-\delta+1)}|\nabla\varphi|^p\Big)\\[5pt]
&\leq &a'\Big(d_0^{\,-p(\alpha-\delta)}|\Delta\chi|^p\,|\eta|^p+d_0^{\,-p(\alpha-\delta+1)}(|\nabla\chi|^p\,|\eta|^p+|\chi|^p\,|\nabla\eta|^p)\Big)
\end{eqnarray*}
with $a, a'>0$.
Therefore one obtains a bound
\begin{equation}
\int_\Omega|R_\alpha|^p\leq a_2\int_0^\infty dr\,r^{p-1}|\chi(r)|^p+a_3\int_0^\infty dr\,r^{p-1}|\chi'(r)|^p+
a_4\int_0^\infty dr\,r^{2p-1}|\chi''(r)|^p
\label{recc4.2}
\end{equation}
with $a_2,a_3,a_4>0$.
This is very similar to the bounds occurring in the proof of Theorem~\ref{tcc4.1} with the exception of the last term which depends 
on $\chi''$.
If this term were absent one could then replace $\chi$ by the sequence of functions  $\chi_n$ used in the proof of the earlier proposition 
to complete the argument that  $\nu_p(\Omega)\leq \left((p-1)(d-k)(d-k+p\,\delta-2p)\,p^{-2}\right)^p$.
But the extra term complicates things. 
In fact the $\chi_n$ used earlier are not even twice  differentiable.
Therefore it is necessary to make a more sophisticated choice.
We now use an argument given in Section~4 of  \cite{RSi3}.

Let $\chi_n$ be the sequence of functions on $\langle0,\infty\rangle$ used in the proof of Theorem~\ref{tcc4.1}.
The derivatives $\chi_n'$ are discontinuous at $n^{-1}$ and at $1$.
The functions $\xi_n=\chi_n^2$ have similar characteristics to the $\chi_n$ except  their derivatives $\xi_n'$ are only discontinuous at $1$.
Therefore we now consider the $\xi_n$ and modify the derivative $\xi'_n$ by the addition of a linear function to remove the discontinuity.
The modifications $\eta_n$ of the derivatives are defined by $\eta_n(s)=0$ if $s\leq n^{-1}$ or $s\geq 1$ and 
\[
\eta_n(s)=\xi_n'(s)-\xi'_n(1)(s-n^{-1})/(1-n^{-1})
\]
if $s\in[n^{-1},1]$.
Now $\eta_n$ is continuous and we set $\zeta_n(s)=\int_0^s\eta_n$ for $s\leq1$ and $\zeta_n(s)=\zeta_n(1)$ if $s\geq1$.
The resulting function $\zeta_n$ is twice-differentiable.
Finally setting $\rho_n=\zeta_n/\zeta_n(1)$ one verifies that $0\leq \rho_n\leq 1$, $\rho_n(s)=0$ if $s\leq n^{-1}$ and $\rho_n(s)=1$
if $s\geq1$.
Moreover, $\lim_{n\to\infty}\rho_n(s)=1$ for all $s>0$.
Finally set $\sigma_n=\rho_n\,\zeta$ where $\zeta$ is the cutoff function  used in  the proof of Theorem~\ref{tcc4.1}.
Now we consider the estimates (\ref{recc4.1}) and (\ref{recc4.2})  with $\chi$ replaced by the sequence $\sigma_n$.

First since $\sigma_n\to 1$ on $\langle0,1]$ as $n\to\infty$ it follows that $\int^\infty_0dr\,r^{-1}|\sigma_n(r)|^p\to\infty$ as $n\to\infty$
but $\int^\infty_0dr\,r^{p-1}|\sigma_n(r)|^p$ is uniformly bounded in $n$.
Moreover, $\sigma_n'=\zeta_n(1)^{-1}\eta_n\,\zeta+\zeta'$ and it follows by the earlier calculation that  $\int^\infty_0dr\,r^{p-1}|\sigma'_n(r)|^p$ is 
also uniformly bounded in $n$.
Therefore it remains to consider the term in (\ref{recc4.2}) dependent on $\sigma_n''$.
But $\sigma''_n=\zeta_n(1)^{-1}(\eta_n'\zeta+\zeta')+\zeta''$.
Therefore it follows from the definition of $\eta_n$  and the cutoff $\zeta$ that 
\[
\int^\infty_0dr\,r^{2p-1}|\sigma_n''|^p\leq a+b\int^1_{n^{-1}}dr\,r^{2p-1}|\xi_n''(r)|^p
\]
with $a,b>0$ independent of $n$.
Now on $[n^{-1},1]$ one has 
\[
\xi''_n(r)=2\,(\chi_n'(r))^2+2\,\chi_n(r)\chi_n''(r)=2\,r^{-2}(1-\log rn)/(\log n)^2
\;.
\]
Therefore
\[
\int^1_{n^{-1}}dr\,r^{2p-1}|\xi_n''(r)|^p=2^p(\log n)^{-2p}\int^1_{n^{-1}}dr\,r^{-1}|1-\log rn|^p\leq 2^{p-1}(\log n)^{-(p-1)}
\]
and this gives a  bound   uniform for $n>1$.

One now deduces that if  $\varphi_\alpha$ in the bound (\ref{recc4.0}) is replaced by $\varphi_{\alpha,n}=d_0^{\,-\alpha}\eta\,\sigma_n$
then in the limit $n\to\infty$ the second term tends to zero since the numerator is bounded uniformly for $n>1$ and the denominator converges to infinity.
Therefore one concludes that 
\[
\nu_p(\Omega)\leq (1-\lambda)^{-(p-1)}\left((p-1)(d-k)(d-k+p\,\delta-2p)\,p^{-2}\right)^p
\] 
for all $\lambda\in\langle0,1\rangle$.
Hence in the limit $\lambda\to0$ one obtains the first bound of the theorem.
This was, however, obtained with the assumption $c(s)=s^\delta$.
But again by rescaling one can arrange that the $\sigma_n$ are supported in a small interval $[0,r]$ and this allows one to reduce the general case to the special case.
There is one extra small complication which did not occur in the Hardy case and that arises since the weighting factor $c_\Omega$ 
is positioned centrally in the operator $H$ and is not a direct weighting of the measure.
But this causes no difficulty. 
For example, if $\varphi$ has support within  distance $r$ of the boundary then
\begin{eqnarray*}
|(Hd_0^{\,-\alpha})\varphi|&\leq&c_\Omega\,| \Delta d_0^{\,-\alpha}|\,|\varphi|
+|c'_\Omega|\,|(\nabla d_0).(\nabla d_0^{\,-\alpha})|\,|\varphi|\\[5pt]
&\leq& \left(d_0^{\,\delta}|\Delta d_0^{\,-\alpha}|
+d_0^{\,\delta-1}|\nabla d_0^{\,-\alpha}|\right)|\varphi|\,(1+r^{|\delta-\delta'|})
\;.
\end{eqnarray*}
Making these modifications one obtains the first bound of the theorem modulo an additional factor $(1+r^{|\delta-\delta'|})$
but since this is valid for all small $r>0$ one can then take the limit $r\to0$.

\medskip

\noindent{\bf Step 2}$\;$ Next we assume $k_\infty=k$ and  establish the second bound in Theorem~\ref{tcc4.31}.
The proof is similar to that of Theorem~\ref{tcc4.3}.

We continue to use the factorization  $\Ri^d=\Ri^k\times \Ri^{d-k}$ and  to  set $x=(y,z)\in \Omega$ with 
$y\in \Ri^k$ and $z\in\Ri^{d-k}$. 
Then  $ d_\Gamma(y,z)=|z|$ if $y\in K$
and  the Rellich  inequality (\ref{ecc1.2}) on $L_p(\Omega)$  takes the form
\begin{equation}
\int_{K}dy\int_{\Ri^{d-k}} dz\,|(H\varphi)(y,z)|^p\geq 
 c_p^{\,p}\int_{K}dy\int_{\Ri^{d-k}} dz\, c(|z|)^p|z|^{-2p}|\varphi(y,z)|^p
\label{ecc4.22}
\end{equation}
for all $\varphi\in C_c^2(K\times \Ri^{d-k})$. 
Therefore the  optimal constant  satisfies
\[
\nu_p(\Omega)\leq
\bigg(\int_{K}dy\int_{\Ri^{d-k}} dz\,|(H\varphi)(y,z)|^p\Big/
\!\int_{K}dy\int_{\Ri^{d-k}} dz\,  c(|z|)^p|z|^{-2p}\, |\varphi(y,z)|^p\bigg)
\]
for all $\varphi\in C_c^2(K\times \Ri^{d-k})$.

Again we set $\varphi=\eta\, \chi$ with $\chi\in C_c^\infty(O_{\!R})$, where $O_{\!R}=\{z\in \Ri^{d-k}:|z|>R\}$, and  $\eta\in C_c^\infty(K)$.
But the action of $H$ on the product $ \chi\,\eta$ takes the  Grushin form
\begin{eqnarray*}
(H\varphi)(y,z)&=&-\sum^k_{j=1}c(|z|)\chi(z)\,(\partial_j^{\,2}\eta)(y)-\sum^d_{j=k+1}(\partial_jc(|z|)\partial_j\chi)(z)\,\eta(y)\\[5pt]
&=&c(|z|) \chi(z)(\Delta\eta)(y)+(H\chi)(z)\eta(y)
\end{eqnarray*}
where the second line is a slight abuse of notation.
This identity replaces the Leibniz rule used in the proof of Theorem~\ref{tcc4.3}.

Then arguing as in the former proof one 
 obtains for all $\lambda\in\langle0,1\rangle$ the estimates
\begin{eqnarray*}
\nu_p(\Omega)&\leq&(1-\lambda)^{-(p-1)}\Bigg({{\int_{O_{\!R}} dz\,|(H\chi)(z)|^p}\over{ \int_{O_{\!R}}dz\,c(|z|)^p|z|^{-2p}|\chi(z)|^p}}\Bigg)
\nonumber\\[5pt]
&&\hspace{3.3cm}+\lambda^{-(p-1)}\Bigg({{\int_{O_{\!R}} dz\,c(|z|)^p|\chi(z)|^p}\over{ \int_{O_{\!R}} dz\,c(|z|)^p|z|^{-2p}|\chi(z)|^p}}\Bigg)
\Bigg({{\int_K dy\,|(\Delta\eta)(y)|^p}\over{ \int_K dy|\eta(y)|^p}}\Bigg)
\end{eqnarray*}
as a replacement for (\ref{ecc4.14}).
But  since $k_\infty=k$ the infimum over $\eta$ of the second term on the  right hand side is zero.
This is no longer a consequence of Lemma~\ref{lcc4.0} but it follows by identical reasoning.
Hence one can then take the limit $\lambda\to0$ to deduce that
\[
\nu_p(\Omega)\leq\Bigg({{\int_{O_{\!R}} dz\,|(H\chi)(z)|^p}\over{ \int_{O_{\!R}}dz\,c(|z|)^p|z|^{-2p}|\chi(z)|^p}}\Bigg)
\]
for all $\chi\in C_c^\infty(O_{\!R})$.
Thus the problem of estimating $\nu_p(\Omega)$ is reduced to a `large distance' estimate on  the Rellich 
constant $\nu_p(\Ri^{d-k}\backslash\{0\})$.
This follows from the standard argument sketched in the proof of Proposition~\ref{pcc4.2}.
One obtains the bound
\[
\nu_p(\Omega)\leq \left((p-1)(d-k)(d-k+p\,\delta'-2p)\,p^{-2}\right)^p
\;.
\]
The second statement of the theorem then follows by minimizing this bound and the local bound obtained in Step~1 of the proof.

\smallskip

The proof of Theorem~\ref{tcc4.31} is completed by  noting that if $\delta=\delta'$ the upper bound on $\nu_p(\Omega)$ coincides with 
the lower bound given by Theorem~\ref{tcc1.2}.
Therefore one has equality between $\nu_p(\Omega)$  and the bound.
\hfill$\Box$

\bigskip

Although Proposition~\ref{pcc4.2} and Theorem~\ref{tcc4.31} do not provide compelling evidence that the optimal constant
in the Rellich inequality should be $C_p^{\,p}$  the arguments of \cite{Rob12} give some support to this conjecture  in more general circumstances.
The following $L_2$-result applies on the complement of a general convex set and for all $\delta, \delta'\in[0,2\rangle$.

\begin{prop}\label{pcc4.3}
Adopt the assumptions of Theorem~$\ref{tcc1.2}$ but with 
 $p=2$. 
 It follows that if $d-d_H+2(\delta\wedge\delta')-4>0$ then  the Rellich inequality $(\ref{ecc1.2})$ is valid with constant equal to  $C_2^{\,2}$.
\end{prop}
\proof\
The proposition is essentially a corollary of Theorem~1.2 in \cite{Rob12}.

First  Theorem~\ref{tcc1.1}  of the current paper establishes the Hardy inequality 
\[
\int_\Omega c_\Omega\,|\nabla\varphi|^2\geq \int_\Omega\,|\eta\,\varphi|^2
\]
with $\eta=a_2\,c_\Omega^{\,1/2}d_\Gamma^{\;-1}$ where $a_2=(d-d_H+(\delta\wedge\delta')-2)/2$.
Secondly
\[
\int_\Omega c_\Omega\,|\nabla\eta|^2=a_2^{\,2}c_\Omega^{\,2}d_\Gamma^{\;-4}\,|1-c'_\Omega d_\Gamma/2c_\Omega|\leq (\nu /a_2^2)\,\eta^4
\]
where $\nu=\sup\{|1-t/2|^2:\delta\wedge\delta'\leq t\leq \delta\vee\delta'\}$.
In particular $\nu=(1-(\delta\wedge\delta')/2)^2$ if $\delta,\delta'\in[0,2\rangle$. 
Theorem~1.2 in \cite{Rob12} asserts, however,  that if $\nu /a_2^2<1$ then the Rellich inequality (\ref{ecc1.2}) is satisfied
with constant $\nu_2(\Omega)=a_2^4(1-\nu /a_2^2)^2=(a_2^2-\nu)^2$.
But  the condition $\nu<a_2^2$ is equivalent to $d-d_H+2(\delta\wedge\delta')-4>0$ or, equivalently, to $C_2>0$. 
Then one calculates that
\[
\nu_2(\Omega)=(a_2^2-\nu)^2=((d-d_H)(d-d_H+2(\delta\wedge\delta')-4)/4)^2
\;.
\]
But the right hand side is equal to $C_2^{\,2}$.
\hfill$\Box$

\begin{remarkn}\label{rcc4.2}  
The proof of the lower bound $\nu_2(\Omega)\geq C_2^{\,2}$ established in Proposition~\ref{pcc4.3}
 readily extends to all $\delta,\delta'\geq0$ with $\delta+\delta'<4$.
Moreover, if $K=\{0\}$ then it also follows from Remark~\ref{rcc4.1}, with $p=2$, that $\nu_2(\Omega)\leq C_2^{\,2}$.
Therefore the conclusion  $\nu_2(\Omega)= C_2^{\,2}$ of Proposition~\ref{pcc4.2} is valid for $\delta\neq\delta'$ if $p=2$.
\end{remarkn}

\section{Concluding remarks}\label{S5}

In conclusion  we note that the $L_2$-Rellich inequalities established in \cite{Rob12} are much stronger than the corresponding $L_2$-statement of  Theorem~\ref{tcc1.2}.
If $p=2$ the Hardy inequality (\ref{ecc1.2})  gives a lower bound on the quadratic form $h(\varphi)=\int_\Omega c_\Omega\,|\nabla\varphi|^2$ which is valid for all $\varphi\in C_c^1(\Omega)$.
But the form $h$ is closeable and the lower bound extends to the closure $\overline h$.
The latter is, however, a local  Dirichlet form and it determines in a canonical manner a submarkovian operator $H_{\!F}$, the Friedrichs' extension
of the symmetric operator $H=-\divv(c_\Omega\nabla)$ defined in the introduction on $C_c^2(\Omega)$.
But the domain of the closed form $\overline h$ is equal to $D(H_{\!F}^{\,1/2})$ and $\overline h(\varphi)=\|H_{\!F}^{\,1/2}\varphi\|_2^2$
for all $\varphi\in D(H_{\!F}^{\,1/2})$.
In particular the $L_2$-Hardy inequality of Theorem~\ref{tcc1.1} can be rephrased in a `weighted operator' form
\[
\|H_{\!F}^{\,1/2}\varphi\|_2^2=\int_\Omega |H_{\!F}^{\,1/2}\varphi|^2\geq a_2^{\,2}\int_\Omega c_\Omega|d_\Omega^{\,-1}\varphi|^2
=a_2^{\,2}\|c_\Omega^{\,1/2}d_\Omega^{\,-1}\varphi\|_2^2
\]
for all $\varphi\in D(H_{\!F}^{\,1/2})$.
It can be stated equivalently as
\[
H_{\!F}\geq a_2^{\,2}\,c_\Omega\,d_\Omega^{\;-2}
\]
 in the sense of ordering of positive self-adjoint operators.
This form of the Hardy inequality is the starting point of Theorem~1.2 of \cite{Rob12}.
The conclusion of the latter theorem is the validity of the Rellich inequality
\begin{equation}
\|H_{\!F}\varphi\|_2^2=\int_\Omega |H_{\!F}\varphi|^2\geq c_2^{\,2}\int_\Omega c_\Omega^2|d_\Omega^{\,-2}\varphi|^2
=a_2^{\,2}\|c_\Omega^{\,1/2}d_\Omega^{\,-1}\varphi\|_2^2
\label{ecc5.1}
\end{equation}
for all $\varphi\in D(H_{\!F})$ or, in the sense of  operator ordering,
 \[
 H_{\!F}^2\geq  c_2^{\,2}\,c_\Omega^{\,2}\,d_\Omega^{\;-4}
 \;.
 \]
But the statement of Theorem~\ref{tcc1.2} in the introduction only gives a
statement comparable to (\ref{ecc5.1}) for $\varphi\in C_c^2(\Omega)$ or, by closure for all $\varphi\in D(\overline H)$ where $\overline H$ is the closure of $H$.
In particular it  gives the operator statement
 \[
 H^*\overline H\geq  c_2^{\,2}\,c_\Omega^{\,2}\,d_\Omega^{\;-4}
 \;.
 \]
Since $H_{\!F}\supseteq \overline H$ it follows that 
 $(H_{\!F})^2\leq  H^*\overline H$ with equality if and only if $H$ is essentially self-adjoint,
 i.e.\ if and only if $H^*=\overline H=H_{\!F}$.
 Hence the $L_2$-Rellich  inequalities of \cite{Rob12} are strictly stronger than those of the current paper
 unless $H$ is essentially self-adjoint.

Another way of distinguishing between the two classes of symmetric operators is to consider the case that 
$H^2$ is densely-defined as a positive symmetric operator.
Then $H^*\overline H$ corresponds to the Friedrichs' extension $(H^2)_{\!F}$ of $H^2$.
For example, if $H=\Delta$, the Laplacian defined on $C_c^2(\Omega)$, then $\Delta^2$ is a symmetric operator
on $C_c^4(\Omega)$.
But  $\Delta_F$ is the the self-adoint extension of $\Delta$ with Dirichlet boundary conditions
and  $\Delta^*\overline \Delta=(\Delta^2)_{\!F}$ is the biharmonic operator which  is determined by 
quite different boundary conditions to those of $(\Delta_F)^2$.
If one considers the classic case  of $\Omega=\Ri^d\backslash\{0\}$ it is well known that $\Delta$ has a unique
submarkovian extension if and only if $d>2$, which happens to be the condition that ensures the validity of the 
Hardy inequality.
Moreover, $\Delta$ has a unique self-adjoint extension if and only if $d>4$, which is the condition which ensures
the validity of the Rellich inequality.
So in this simple case there is no ambiguity. 
But for more general $\Omega$ and operators $H=-\divv(c_\Omega\nabla)$ the situation is much more complicated.
Criteria for Markov uniqueness were obtained for quite general $\Omega$ in \cite{LR} and it would be of interest 
to develop analogous criteria for self-adjointness.

\section*{Acknowledgements}
The author is again  indebted to Juha Lehrb\"ack for  continuing advice and information on  Hardy and Rellich inequalities
and for a critical reading of a  preliminary version of the paper.


\end{document}